\documentclass[12pt]{article}

\usepackage[table,xcdraw]{xcolor}
\usepackage[utf8]{inputenc} 
\usepackage[T1]{fontenc}    

\usepackage{url}            
\usepackage{booktabs}       
\usepackage{amsfonts}       
\usepackage{nicefrac}       
\usepackage{microtype}      
\usepackage{lipsum}
\usepackage{graphicx}
\usepackage{bbm}

\usepackage{caption} 
\usepackage{natbib}

\usepackage{amssymb}
\usepackage{amsmath}
\usepackage{mathtools}
\usepackage{amsthm}
\usepackage{multirow}
\usepackage{hyperref}

\hypersetup{   
    colorlinks=true,
    linkcolor=blue,
    filecolor=magenta,      
    urlcolor=cyan,
    citecolor = blue,
} 

\usepackage[hmarginratio=1:1, vmarginratio=1:1, marginparwidth=4cm, textheight=23cm, textwidth=18cm]{geometry}

\def \aic{\mathrm{AIC}}
\def\bbp{1 + \sqrt{c}}
\newtheorem{theorem}{Theorem}[section]
\newtheorem{lemma}[theorem]{Lemma}

\title{\textbf{High dimensional PCA: a new model selection criterion}}

\author{
 \normalsize{\textbf{Abhinav Chakraborty}}\\
 \normalsize{Department of Statistics}\\
 \normalsize{The Wharton School}\\
 \normalsize{University of Pennsylvania}\\
 \normalsize{Philadelphia, PA 19104, USA}\\
 \normalsize{\texttt{abch@wharton.upenn.edu}}
 
 \and
 
 \normalsize{\textbf{Soumendu Sundar Mukherjee}}\\
 \normalsize{Interdisciplinary Statistical Research Unit}\\
 \normalsize{Applied Statistics Division}\\
 \normalsize{Indian Statistical Institute}\\
 \normalsize{Kolkata, WB-700108, India}\\
 \normalsize{\texttt{soumendu041@gmail.com}}
 
 \and
 
 \normalsize{\textbf{Arijit Chakrabarti}}\\
 \normalsize{Applied Statistics Unit}\\
 \normalsize{Applied Statistics Division}\\
 \normalsize{Indian Statistical Institute}\\
 \normalsize{Kolkata, WB-700108, India}\\
 \normalsize{\texttt{arc@isical.ac.in}}
}

\usepackage{fancyhdr}
\begin{document}
\maketitle
\begin{abstract}
Suppose we have a random sample from a multivariate population consisting of many variables. Estimating the number of dominant/large eigenvalues of the population covariance matrix based on the sample information is an important question arising in Statistics with wide applications in many areas. In the context of Principal Components Analysis (PCA), the linear combinations of the original variables having the largest amounts of variation are determined by this number. In this paper, we study the high dimensional asymptotic setting where the number of variables grows at the same rate as the number of observations. We work in the framework where the population covariance matrix is assumed to have the spiked structure proposed in \cite{johnstone2001distribution}). In this setup, the problem of interest becomes essentially one of model selection and has attracted a lot of interest from researchers. Our focus is on the Akaike Information Criterion (AIC) which is known to be strongly consistent from the work of \cite{bai2018consistency}. The result of \cite{bai2018consistency} requires a certain ``gap condition'' ensuring that the dominant eigenvalues of the covariance matrix are all above a level which is strictly larger than a threshold discovered by Baik, Ben Arous and Peche (called the BBP threshold), both quantities depending on the limiting ratio of the number of variables and observations.
It is well-known in the literature that, below this threshold, a spiked covariance structure becomes indistinguishable from one with no spikes. Thus the strong consistency of AIC requires in a sense some extra ``signal strength'' than what the BBP threshold corresponds to. \\

In this paper, our aim is to investigate whether consistency continues to hold even if the ``gap'' is made smaller. In this regard, we make two novel theoretical contributions. Firstly, we show that strong consistency under arbitrarily small gap is achievable if we alter the penalty term of AIC suitably depending on the target gap. Inspired by this result, we are able to show that a further intuitive alteration of the penalty can indeed make the gap exactly zero, although we can only achieve weak consistency in this case. We compare the two newly-proposed estimators with other existing estimators in the literature via extensive simulation studies, and show, by suitably calibrating our proposals, that a significant improvement in terms of mean-squared error is achievable.
\end{abstract}


\section{Introduction}
 Suppose we have a sample of observations from a multivariate population with $p$ variables. Estimating the number of dominant/significant eigenvalues of the population covariance matrix in such a scenario is an important question arising in Statistics with wide applications in many areas. In the context of Principal Components Analysis (PCA), a very popular method of dimension reduction for multivariate data, the individual principal components having the largest variability are determined by this number.

Let $\lambda_1 \geq \lambda_2 \geq \cdots \geq \lambda_p$ be the eigenvalues of the population covariance matrix. We further assume that the covariance matrix has a spiked structure proposed by \cite{johnstone2001distribution}. In this framework, the number of dominant eigenvalues is denoted by $k$ and all the eigenvalues except the first $k$ are all assumed equal, i.e. $\lambda_1 \geq \lambda_2 \geq \cdots \geq\lambda_k > \lambda_{k+1} = \lambda_{k+2} = \cdots = \lambda_{p}$. This $k$ is called the true number of dominant/significant components in this framework.

The spiked covariance model  finds wide applications in many scientific fields. In wireless communications, for example, a signal emitted by a source is modulated and received by several antennas, and the  quality of reconstruction of the original signal is directly linked to the ``inference'' of spikes. The spiked model  is also used in different areas of Artificial Intelligence such as face, handwriting and speech recognition and statistical learning. See  \cite{johnstone2018pca} for more applications.  

The number of significant components $k$ is usually unknown, and we need to estimate it,  which in turn becomes essentially a problem of model selection in the spiked covariance framework. This will be explained clearly in the next section. Many estimators have been developed  in the literature, mostly based on information theoretic criteria, such as the minimum description length (MDL), Bayesian Information Criterion (BIC) and Akaike Information Criterion (AIC) (see, e.g., \cite{wax1985detection}).  However, these applications have been focused on the large sample size and low dimensional regimes and arguments in support of these estimators may not carry over to the high dimensional setup. In the recent past, several works have appeared in the area of signal processing for high dimensional data, where techniques from random matrix theory (RMT) it have been used (see, for example, \cite{kritchman2009non} and \cite{nadler2010nonparametric}). More recent papers in the literature include  \cite{passemier2014estimation} and \cite{bai2018consistency}.        

\begin{figure}[!htbp]
\centering
    \includegraphics[scale =0.6]{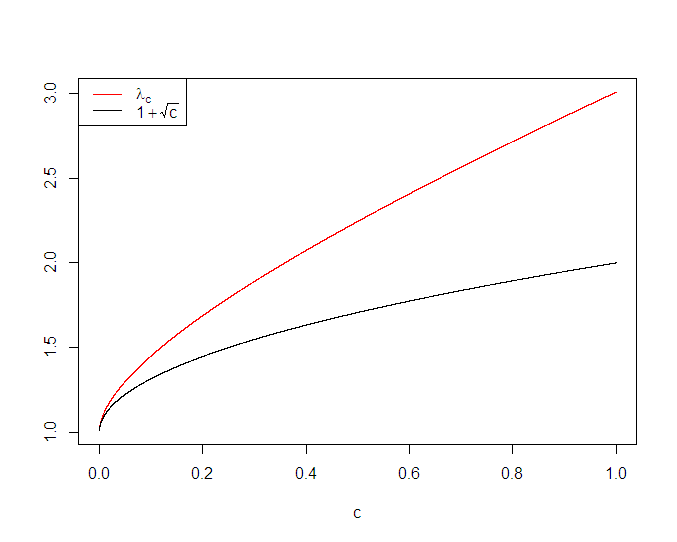}
    \caption{Gap between $\lambda_c$ and the BBP threshold.}
    \label{fig:gap}
\end{figure}

 Our work is inspired by \cite{bai2018consistency}, where the authors consider the Akaike Information Criterion (AIC) (\cite{akaike1998information}) and the Bayesian Information Criterion (BIC) (\cite{schwarz1978estimating}) as their estimation criterion in a high dimensional setting. 
They studied the consistency of the estimators based on the AIC and BIC criteria under an asymptotic framework where $p,n \to \infty$ such that $p/n \to c >0$,  where $n$ refers to the sample size. They showed that unlike AIC, the BIC criterion is not consistent when the signal strength $\lambda_k$ is bounded (\cite{bai2018consistency});
in other words BIC requires much more signal strength for signal detection. Their main result shows consistency of AIC for estimating $k$, when a certain ``gap condition'' is satisfied, i.e. when $\lambda_k$ is above a certain level $\lambda_c$ above the \textit{BBP threshold }(\cite{baik2005phase}) of $\bbp$ (see  Section \ref{sec:baistrongconsistency}). More precisely, they showed that AIC is consistent if and only if $\lambda_k >\lambda_c>1+\sqrt{c}$. Figure \ref{fig:gap} shows the ``gap'' between $\lambda_c$ and $1+\sqrt{c}$. We want to highlight the fact that if $\lambda_k \leq 1+\sqrt{c}$, then there is no hope for estimating $k$, (see, for example, \cite{baik2005phase} and  Section \ref{sec:rmt_results} for more details).

The primary aim our work is to investigate if we can improve upon the results of \cite{bai2018consistency} by suitably modifying the AIC criterion to give consistent estimates of $k$ under weakening of the ``gap condition''. Towards this, our main contributions are as follows. We have shown that, given any $\delta>0$, we can develop an estimator (depending on $\delta$) which is strongly consistent when the gap between $\lambda_k$ and the BBP threshold is at least $\delta$ (see Section \ref{sec:strongconsistencyresult}). We can also make the gap exactly zero by modifying our estimator, but then we have to part with strong consistency-we can only prove weak consistency of the modified estimator (see Section \ref{sec:weakconsistencyresult}).  We note that there is another weakly-consistent estimator due to \cite{passemier2014estimation} that is known in the literature and which also works under zero gap.

The inspection of the proof of strong consistency in \cite{bai2018consistency} reveals that for any arbitrary $\delta > 0$, a modification (based on $\delta$ and ratio of dimension and sample size) of the penalty of AIC will make the asymptotic argument work when the gap between $\lambda_k$ and the BBP threshold is above $\delta$. A further modification of the penalty term by letting that $\delta$ go to zero at appropriate rate depending on $n$ gives us an estimator that is weakly consistency under `zero gap'. This proof is obtained by employing novel arguments aided by some deep Random Matrix Theory results on the asymptotic behavior of sample covariance matrices.

After finishing this work, it came to our notice that a similar problem has been studied independently by \cite{hu2020detection}. However, their motivations, results and methods of proof are different from ours.

The rest of the paper is organised as follows. In Section \ref{sec:litreview}, we discuss the problem setup, review some key results of Random Matrix Theory (RMT) and the literature on estimating the number of spikes. Section \ref{sec:ideasResult} introduces the main idea and results behind our modified AIC criterion. In Section \ref{sec:simulation} we compare our proposed estimator with other estimators available in the literature via extensive simulation studies. Conclusions then follow and the Appendix collects all the proofs.

\section{Problem setup and a review of relevant literature}
\label{sec:litreview}
 We first describe our setup and main problem in the first few paragraphs of this section. We follow the same notations as used in \cite{bai2018consistency}. Suppose we have a random sample $\mathbf{y_1,\hdots,y_n}$ of size $n$ from a population of dimension $p$ and let $\mathbf{Y} = (\mathbf{y_1,\hdots,y_n})^{T}$ denote the full $(n \times p)$ data matrix. Let the population mean be $\mathbf{\mu}$ and population covariance matrix is denoted by $\mathbf{\Sigma}$. Our interest here is in the spiked covariance structure described before which can be succinctly described through models $M_{k}$, where  
\[
M_k : \lambda_1 \geq \lambda_2 \geq \cdots \geq \lambda_k > \lambda_{k+1} = \cdots = \lambda_p. 
\]
We consider models $M_{0}, M_{1}, \ldots, M_{p-1}$ starting from a configuration of zero spikes (i.e. $\lambda_1=\cdots=\lambda_p$) to one with $(p-1)$ spikes. We want to estimate, using the available data, the true value of $k$. So our main problem is thus reduced to one of model selection from the pool of the above $p$ candidate models.

 Our focus in this paper will be the high dimensional asymptotic setting where $p$ and $n$ grow proportionately, i.e. we assume that
\begin{equation}
    p/n \to c >0. \tag{C1}\label{eq:C1}
\end{equation}
To avoid notational clutter, we suppress the dependence of $p$ on $n$, and of $\mathbf{\Sigma}$,$\lambda_1,\lambda_2 \hdots $ on $p$. For the estimation of $k$, we will be interested in the study (to be explained later in this section) of the distributional properties of the eigenvalues of the sample covariance matrix of the $\mathbf{y_j}$’. We can assume without any loss of generality that $\mathbf{\mu} =0$. Then the sample covariance matrix $\mathbf{S_n}$ is given by  
\[
\mathbf{S_n} = \frac{1}{n-1}\sum_{i=1}^n\mathbf{y_i} \mathbf{y_i}^T. 
\] 
Denoting $\{1,\hdots,n\}$ by $[n]$, we may also assume, as in \cite{bai2018consistency}, that $\mathbf{y_j}=\mathbf{\Sigma}^{\frac{1}{2}} \mathbf{x_j}$ for $j\in [n]$, where $\mathbf{x_k}=(x_{1k},\hdots,x_{pk})^T$, and $\{x_{ij},i \in [n],j \in[n]\}$ is a double array of i.i.d. random variables with mean $0$, variance $1$. We further assume, as in \cite{bai2018consistency}, that the $x_{ij}$'s have finite  with finite fourth moment. 

For the rest of the paper we assume that the eigenvalues of population covariance matrix $\mathbf{\Sigma}$ are 1 except the first $k$ which are $(\lambda_i)_{1\leq i\leq k}$ and $\mathbf{\Sigma}$ has the form
\begin{equation}
        \mathbf{\Sigma}  =
\begin{pmatrix}
\mathbf{\Sigma_k} & \mathbf{0}\\
\mathbf{0} & \mathbf{I_{p-k}}
\end{pmatrix} \tag{C2}\label{eq:C2}
\end{equation}

where $\mathbf{\Sigma_k}$ has $k$ non-null and non-unit eigenvalues $\lambda_1\geq\lambda_2\geq\cdots\geq\lambda_k>1$. We also assume, as in \cite{bai2018consistency} that $k$ is unknown but a fixed finite number and does not change with $n$.

We will now summarize in the next few subsections, the concepts, issues and results from the literature which are relevant in this study. We will first touch upon some fundamental results on the properties of the eigenstructure of the sample covariance matrix and the high dimensional phenomenon. This will be followed by review of some specific works in the literature on the topic of our interest.

\subsection{Basic results on the eigen-structure of sample covariance matrix}

We start this section by briefly describing Principal Components Analysis (PCA) since the eigenstructure of the sample covariance matrix is intrinsically related to this. We recall that in PCA (see, e.g., \cite{anderson1958introduction, jolliffe1986principal}), one sequentially finds orthogonal directions to produce uncorrelated (normalized) linear combinations of the original variables with maximum variance. In the standard approach to this, this is obtained from the eigen-decomposition of the (sample)covariance matrix. Consider the eigen-decomposition of $\Sigma$:
$$\mathbf{\Sigma} = \lambda_1 \mathbf{u_1 u_1'} + \cdots +\lambda_p \mathbf{u_p u_p'} = \mathbf{U\Lambda U'},$$
where $\mathbf{U}$ is a $p \times p$ orthogonal matrix whose columns are the eigenvector $\mathbf{u_i}$ and $\mathbf{\Lambda}$ is a diagonal matrix with entries $\lambda_i$ being the eigenvalues of $\mathbf{\Sigma}$. The sample analogue of this is
$$\mathbf{S_n} = l_1 \mathbf{v_1 v_1'} + \cdots +l_p \mathbf{v_p v_p'} = \mathbf{VLV}',$$
where now the orthogonal matrix $\mathbf{V}$ has columns which are the sample eigenvectors $\mathbf{v_i}$ and $\mathbf{L}$ is a diagonal matrix consisting of the eigenvalues $l_i$ of $\mathbf{S_n}$. If $\mathbf{Z}$ denotes a generic random observation from the distribution, then the vector of the population principal components is given by $\mathbf{ZU}$. The sample principal components are defined as $\mathbf{Y v_i}$, for $i=1,\cdots,p$. It may be noted that for each $i \in \{1,\cdots,p\}$, the variance of $\mathbf{Z u_i}$ is $\lambda_i$. So larger the $\lambda_i$, larger is the variance of random variable $\mathbf{Z u_i}$. Thus the number of dominant/large eigenvalues of $\mathbf{\Sigma}$ correspond to the number of principal components having the largest amount of information. In the traditional fixed dimensional setting where $p$ is fixed, and $n$ is large, the sample eigenvalues and eigenvectors converge to their population counterparts, i.e. as $n \to \infty$, we have that $l_i \overset{a.s.}{\to} \lambda_i$ and $v_i \overset{a.s.}{\to} u_i$, for each $i = 1,\cdots,p$ (see, for example,  \cite{anderson1958introduction}). The situation is more subtle in the high dimensional case where $p$ grows with $n$ and it is well-known in the literature (\cite{silverstein1995strong})  that the consistency of the sample eigenvalues and eigenvectors do not carry over in this case. In the next part of this section we focus on the high-dimensional aspect with special emphasis on results from random matrix theory under the spiked covariance framework.

\subsubsection{The high dimensional phenomenon and Random Matrix Theory results}\label{sec:rmt_results}
Many of the things described below are classically well-known facts from Random Matrix Theory (RMT) and some are more recent interesting discoveries. The notational convention is similar to \cite{bai2018consistency}. Recall that the eigenvalues of $\mathbf{S_n}$ are denoted by $l_{1}\geq l_{2}\geq \cdots \geq l_{p} \geq 0$ ( again suppressing their dependence on $n$ and/or $p$). Let us define the Empirical Spectral Distribution (ESD) of $\mathbf{S_n}$ by
$$F_n(x) = \frac{1}{p}\sum_{i=1}^p \mathbbm{1}\{l_{i} \leq x\}.$$
By a result of \cite{silverstein1995strong} we have that under condition (\ref{eq:C2}) the ESD of $\mathbf{S_n}$, i.e. $F_n(x)$ converges to $F_{c}(x)$ almost surely, where $F_{c}(x)$ is the Mar\v chenco-Pastur law/distribution. 
Here, for $0<c\leq 1$, $F_{c}$ is given by
$$F_c'(x) = f_c(x) = \frac{1}{2 \pi x c}\sqrt{(b-x)(x-a)}.$$
The support of this distribution is  $[a,b]$ where $a := (1-\sqrt{c})^2$ and $b := (1+\sqrt{c})^2$.

If $c>1$, $F_c$ has a point mass $1 - 1/c$ at the origin, i.e.
$$
F_c(x)=
\begin{cases}
0  &\text{if } x<0,\\
1-1/c &\text{if } 0\leq x<a,\\
1-1/c+\int_a^x f_c(t)dt &\text{if } a\leq x\leq b,
\end{cases}
$$
where $a$ and $b$ are the same as in the case $0<c\leq 1$. 

The previous result characterizes the bulk behaviour of the sample eigenvalues. Now we state results characterizing convergence of individual eigenvalues of $\mathbf{S_n}$. {\it As in \cite{bai2018consistency}, we define $\lambda_i$ to be a \emph{``distant spiked eigenvalue''} if $\lambda_i > 1+\sqrt{c}$.} We also define the function $\psi_c(x)$ for $x \neq 1$ as
$$\psi_c(x) =x + \frac{cx}{x-1}.$$
This next result is the same as Lemma 2.1 in \cite{bai2018consistency}.
\begin{lemma}
\label{asconvergence}
Let $l_{i}$ denote the i-th largest eigenvalue of $\mathbf{S_n}$ the covariance matrix in our setup. Suppose that $\mathbb{E}(x_{11}^4) <\infty $, conditions $(\ref{eq:C1})$ and $(\ref{eq:C2})$ hold, and that $\lambda_1$ is bounded.
\begin{itemize}
    \item[(i)] If $\lambda_i$ is distant-spiked, then $l_{i} \overset{a.s.}{\to}\psi_c(\lambda_i) = \psi_i = \lambda_i + \frac{c\lambda_i}{\lambda_i -1}$.
    \item[(ii)] If $\lambda_i$ is not distant-spiked and $i/p \to \alpha$, then $l_{i}\overset{a.s.}{\to} \mu_{1-\alpha}$, where $\mu_{\alpha}$ is the $\alpha$-quantile of the MP distribution. In particular, if $i=o(p)$, then $l_{i}\overset{a.s.}{\to} \mu_1 =b= (1+\sqrt{c})^2$.  
\end{itemize}
\end{lemma}

The results in the above lemma are examples of a general high dimensional ``phase transition'' phenomenon observed by many authors (e.g., \cite{bai2012sample}, \cite{baik2005phase}), and often referred as the BBP phase transition phenomenon after \cite{baik2005phase}.  In a nutshell, as summarized in \cite{paul2007asymptotics}, this refers to the fact that if the non-unit eigenvalues of a spiked model are close to one, then their sample counterparts would asymptotically behave as if the true covariance matrix were the identity matrix. However, the asymptotics would change critically if the dominant eigenvalues are larger than the threshold of $(1+\sqrt{c})$. We describe its details and the implications on our assumptions. \\
For understanding its effect on distributional convergence, we assume now for simplicity that $\mathbf{\Sigma}$ is diagonal with a single spike ($k=1$), so that $\mathbf{\Sigma}=\mathrm{diag}\{\lambda_1,1,\cdots,1\}$. When $\lambda_1 = 1$, the largest sample eigenvalue is located near the upper edge $b$ of the MP distribution and fluctuates on the (small) scale $n^{-2/3}$ approximately according to the real valued Tracy-Widom distribution: 
\[
    n^{\frac{2}{3}}\frac{l_1-\mu(c)}{\sigma(c)} \overset{d}{\to} TW_1,
\]
where $\mu(c) = b$ and $\sigma(c) = (\bbp)^{4/3}c^{-1/6}$ and $TW_1$ is a random variable following real valued TW distribution (see \cite{johnstone2001distribution} for more details). For $\lambda_1 \leq \bbp$, the largest sample eigenvalue has the same limiting Tracy-Widom distribution-the small spike in the top population eigenvalue has no limiting effect on the distribution  of sample top eigenvalue. Put in another way, asymptotically the largest sample eigenvalue is of no use in detecting a subcritical spike in the largest population eigenvalue. A phase transition occurs at $\bbp$ : for larger values of $\lambda_1$, the largest sample eigenvalue $l_1$ now has a limiting
Gaussian distribution (See \cite{paul2007asymptotics}), with scale on the usual order of $n^{-1/2}$. The mean of this Gaussian distribution (=$\psi(\lambda_1)$) shows a significant upward bias, being significantly larger than the true value of $\lambda_1$. We now come back briefly to the manifestation of the phase transition in the issue of pointwise convergence as in Lemma \ref{asconvergence}, where we consider a spiked model of $k$ spikes. For the non spiked eigenvalue $\lambda_{k+1}$, the sample counterpart $l_{k+1}$ converges a.s. to $(\bbp)^2$ (by the Lemma \ref{asconvergence} (ii)). If $1<\lambda_k \leq \bbp $ then $l_{k}$ also converges a.s. to $(\bbp)^2$. So asymptotically it becomes difficult to distinguish $l_{k}$ and $l_{k+1}$ (i.e. model $M_k$ and $M_{k+1}$ respectively). This is not the case if $\lambda_k > 1+\sqrt{c}$ as seen in Lemma \ref{asconvergence} (i)). Keeping in mind this phase transition behaviour, we have assumed for our results presented later that the first $k$ eigenvalues are ``distant spiked eigenvalues''

\subsection{Strongly consistent estimators based on AIC}
\label{sec:baistrongconsistency}
This section is mainly based on the paper by \cite{bai2018consistency}. They consider the spiked covariance structure mentioned before which is expressed by models denoted by $M_k$'s for varying $k$ depending on the value of the true number of spikes $k$. To estimate $k$, they propose using the traditional AIC criterion to select a model from among the pool of candidate models and thereby obtain an estimator which is strongly consistent. We shall discuss the $p\leq n$ case first. Defining
\[
C_{p,n} = n\log((n-1)/n)^p + np\{1+\log(2\pi)\},
\]
it has been noted in in \cite{fujikoshi2011multivariate} that the criterion value under model $M_j$ is given by 
\[
\aic_j = n \log (l_{1}\hdots l_{j}) + n(p-j) \log \bar{l}_{j} + 2d_j + C_{p,n},
\]
where  $l_{1}>\hdots >l_{p}$ are the sample eigenvalues of $\mathbf{S_n}$ and for $1\leq j \leq p-1$, $\bar{l}_{j}$ is the arithmetic mean of $l_{j+1},\hdots, l_{p}$, that is,
\[
\bar{l}_{j} = \frac{1}{p-j}\sum_{t = j+1}^p l_{t}.  
\]
Furthermore, $d_j$ denotes the number of independent model parameters under model $M_j$ and is given by 
\begin{align*}
d_j &= pj - \frac{1}{2}j(j+1)+j+1+p\\
&=(j+1)(p+1-j/2).
\end{align*}
The expression of $d_j$ is obtained by looking at the eigen-decomposition of the covariance  matrix $\Sigma$ which can be written as:
\[
\Sigma = \sum_{i=1}^j \lambda_i u_i u_i^T + \lambda(I - \sum_{i=1}^j u_i u_i^T), 
\]
where $u_i$'s are mutually orthogonal unit vectors. It is evident that $p$ degrees of freedom are accounted for by $\mu$, $j+1$ by $\lambda_i$'s $i=1,\hdots j$ and $\lambda$, and $pj - j(j+1)/2$ by the $j$ orthonormal eigen-vectors. The AIC criterion selects the model $M_{\hat{k}_A}$ where
\[
\hat{k}_A = \text{arg min}_j \aic_j.
\]
The estimator of $k$ proposed by \cite{bai2018consistency} is $ \hat{k}_{A}$. When we are interested in only the first $q$ models $M_j, j=0,1, \hdots, q-1$, then the criteria is defined by considering the minimum only with respect to $j=0,1,\hdots,q-1$. We call $q$ the number of candidate models.  Denoting $A_j =\frac{1}{n} (\aic_j - \aic_{p-1})$, the model selection rule of AIC can be equivalently obtained as 
$$
\hat{k}_A = \text{arg min}_j A_j.
$$
Note that $A_j$ is given by
$$
A_j = (p-j)\log \bar{l}_{j} - \sum_{i=j+1}^p \log l_{i} - \frac{(p-j-1)(p-j+2)}{n}.
$$
We are now going to state main results of \cite{bai2018consistency} regarding consistency of $\hat{k}_A$. A criterion $\hat{k}$ for estimating $k$ is said to be consistent (strongly consistent) if $\lim_{n\to\infty} P(\hat{k} = k ) =1$ [ $P(\lim_{n\to\infty}\hat{k} = k ) =1$.

We will first state the result for the case $0<c \leq 1$. Before stating the result, we state condition (\ref{eq:C3}) referred to as the ``gap condition''. Recalling the function $\psi_c(x)$ and quantities $\psi_i=\psi_c(\lambda_i)$ defined earlier, the ``gap condition'' is given by 
\begin{equation}
    \psi_k - 1 - \log \psi_k - 2c > 0.  \tag{C3}\label{eq:C3}
\end{equation}
The next result is proved as Theorem 3.1(i) of \cite{bai2018consistency}.
\begin{theorem}
Suppose the conditions (\ref{eq:C1}) with $0<c \leq 1$, and (\ref{eq:C2}) hold, and that the number of candidate models, $q$, satisfies $q=o(p)$. Suppose also that $\lambda_1$ is bounded.We have the following results on the consistency of the estimation criterion $\hat{k}_A$
based on AIC:\\
 
\begin{itemize}
    \item[(i)] If the gap condition (\ref{eq:C3}) does not hold, then $\hat{k}_A$ is not consistent.
    \item [(ii)] If the gap condition (\ref{eq:C3}) holds, then $\hat{k}_A$ is strongly consistent.
\end{itemize}
\end{theorem}

 Next  we consider the case where $p,n \to \infty$, such that $p>n$ and $p/n \to c >1$. Clearly in this setup the smallest $p -(n-1)$ eigenvalues of $\mathbf{S_n}$ are zero, that is,
\[
l_{n-1} > l_{n} = \hdots = l_{p} = 0.
\]
Thus, as noted in \cite{bai2018consistency},  it is impossible to infer about the smallest population eigenvalues $\lambda_n,\lambda_{n+1},\hdots ,\lambda_p >0$ using the sample eigenvalues. Therefore, these authors additionally assume in this setup that (\ref{eq:C4}) holds where
\begin{equation}
        \lambda_{n-1} = \lambda_{n} = \hdots = \lambda_p =1. \tag{C4}\label{eq:C4}
\end{equation}

 Under this new assumption, for $j=0,1,\hdots,n-2$, it follows that  
 \[
 \tilde{M}_j : \lambda_j > \lambda_{j+1} = \hdots = \lambda_{n-1} \Leftrightarrow M_j: \lambda_j > \lambda_{j+1} = \hdots =\lambda_p. 
 \]
 In order to describe the model selection criterion proposed in this setup in \cite{bai2018consistency}, $\bar{l}_{j}$ first has to be redefined as below:
\[
\bar{l}_{j} = \frac{1}{n-1-j}\sum_{t=j+1}^{n-1} l_{t}   \quad \quad  j= 1,2,\hdots,n-1.
\]
Using this modification, a new criterion $\tilde{A_j}$ is then introduced as
$$
\tilde{A_j} = (n-1-j)\log \bar{l}_{j} - \sum_{i=j+1}^{n-1}\log l_{i} - \frac{(n-j-2)(n-j+1)}{p},
$$
by replacing the $p$ and $n$ in $A_j$ by $n-1$ and $p$, respectively. Note that $\tilde{A}_{n-2} = 0$. The ``quasi-AIC'' rule, henceforth abbreviated as the qAIC rule, selects the model $\hat{k}_{\tilde{A}}$ defined as
\[
\hat{k}_{\tilde{A}} = \text{arg min}_{j\leq n-2} \tilde{A}_j.
\]
The strong consistency result of qAIC for the case $c>1$ is proved under the modified ``gap condition''
\begin{equation}
    \psi_k/c - 1 - \log (\psi_k/c) - 2/c > 0.  \tag{C5}\label{eq:C5}
\end{equation}
We now state the result regarding consistency of $\hat{k}_{\tilde{A}}$ proved as Theorem 3.3(i) in \cite{bai2018consistency}.
\begin{theorem}
Suppose the conditions (\ref{eq:C1}) with $c>1$, and (\ref{eq:C4}) hold, and that the number of candidate models, $q$, satisfies $q=o(p)$. Suppose also that $\lambda_1$ is bounded. We have the following results on the consistency of the estimation criterion $\hat{k}_{\tilde{A}}$
based on qAIC:\\

\begin{itemize}
    \item[(i)] If the gap condition (\ref{eq:C5}) fails, then $\hat{k}_{\tilde{A}}$ is not consistent.
    \item[(ii)] If the gap condition (\ref{eq:C5}) holds, then $\hat{k}_{\tilde{A}}$ is strongly consistent.
\end{itemize}
\end{theorem}

Next we will discuss another type of estimators available in the literature which are weakly consistent, i.e. if $\hat{k}$ is an estimator of $k$ then $\hat{k}\overset{p}{\to} k$ as $n \to \infty$. 
\subsection{Weak consistency}
\label{sec:Passemier and Yao}
 This section is based on the work of  \cite{passemier2014estimation} where they proposed a weakly consistent estimator of $k$ under the ``zero gap'' condition. Suppose we have observed a random sample $\mathbf{y_i, \ldots,y_n}$ from a $p$-dimensional population. These authors additionally assumed that the $\mathbf{y}$'s can be expressed as $\mathbf{y} = \mathbf{EV}^{\frac{1}{2}} \mathbf{x}$, where  $ \mathbf{x} \in \mathbb{R}^p$ is a zero-mean random vector of i.i.d. components, $\mathbf{E}$ is an orthogonal matrix and
\[
\mathbf{V} = \mathrm{cov}(\mathbf{x}) =
\begin{pmatrix}
    \mathbf{\Sigma_k} & \mathbf{0}\\
    \mathbf{0} & \mathbf{I_{p-k}}
\end{pmatrix},
\]  
where $\mathbf{\Sigma_k}$ has $k$ non-null and non-unit eigenvalues $\lambda_1\geq\lambda_2\geq\cdots\geq\lambda_k>1$.
The sample covariance matrix is taken to be 
$$\mathbf{S_n} = \frac{1}{n}\sum_{i=1}^n \mathbf{y_i y_i}'.$$
The proposed estimator in \cite{passemier2014estimation} is based on the differences between consecutive eigenvalues of $\mathbf{S_n}$. The main idea behind that is explained next. Define $\delta_{j}$ as
$$\delta_{j} = l_{j} - l_{j+1}, \mbox{ for }j=1,\cdots,p-1.$$
The authors then observe that under certain assumptions it can be shown that if $j\geq k$ then $\delta_{j} \to 0$ at a fast rate, whereas if $j<k$ and $\lambda_j$'s are different then $\delta_{j}$ tends to a positive limit. Even if $\lambda_j$'s are same, for $j< k$, the convergence of $\delta_{j}$ to zero is slow. Thus a possible estimate of $k$ can be the index $j$ where $\delta_{j}$ becomes small for the first time. The estimator of $k$ proposed in \cite{passemier2014estimation} is denoted by $\hat{k}_P$ and is given by 
$$
\hat{k}_P = \min \{j \in \{1,\hdots,s\}: \delta_{j+1,p} < d_n\},
$$
where $s > k$ is a fixed number big enough, and $d_n$
is an appropriately chosen small number. In practice, the integer
$s$ should be thought as a preliminary bound on the number of
possible spikes. Before stating the main theorem on weak consistency of this estimator, we state one of the main assumptions which is required to prove their theorem.
\vskip5pt
\noindent
\textbf{Assumption \hypertarget{assump1}{1}:} The entries $x^i$ of the random vector $\mathbf{x}$ have a symmetric law and a sub-exponential decay, that is there exists positive constants $C, C'$ such that, for all $t\geq C'$
$$\mathbb{P}(|x^i| \geq t^C) \leq e^{-t}. $$
\begin{theorem}
Let $(\mathbf{y_i})_{(1\leq i\leq n)}$ be $n$ i.i.d. of $\mathbf{y} = \mathbf{EV}^{\frac{1}{2}}\mathbf{x}$, where $\mathbf{x}\in \mathbb{R}^p$ is a zero-mean random vector of i.i.d. components which satisfies  $\hyperlink{assump1}{\text{Assumption 1}}$. Assume that $V$ is of the form described before
where $\mathbf{\Sigma_k} $ has $k$ non null, non unit eigenvalues satisfying
$\lambda_1\geq \hdots\geq\lambda_k>1+\sqrt{c}$. Assume further that (C1) holds.
Let $d_n$ be a real sequence such that $d_n = o(n^{-\frac{1}{2}})$ and  $n^{2/3}d_n\to \infty$, then the estimator $\hat{k}_P$ is a weakly consistent estimator, i.e. $\mathbb{P}(\hat{k}_P = k)\to 1$ as $n \to \infty$. 
\end{theorem}

\section{Main ideas and results}
\label{sec:ideasResult}

In this section we will motivate and describe the main results of this article. We will first discuss our work on strong consistency followed by that on weak consistency. 

\subsection{Strong consistency for modified AIC criterion}
\label{sec:strongconsistencyresult}
The idea for developing a modification of the AIC criterion comes from a careful study of the proof of strong consistency of the traditional AIC in \cite{bai2018consistency}. The ``gap condition'' mentioned before is a crucial requirement in their proof. This makes AIC inadequate for consistent estimation of $k$ when $\lambda_k$ is close to the BBP threshold ($\bbp$). We show in this paper that by suitably modifying the penalty term in the AIC criterion, this problem can be taken care of, in that the modified criterion is strong consistent even when $\lambda_k$ is arbitrarily close to $\bbp$. In this section We  will informally sketch how to get this modified criterion and also the proof of its strong consistency. Towards that, let us first recall the functions we defined earlier 
\begin{align*}
\label{eq:*}
    h(x) &= x - 1 - \log x, \\
    \psi_c(x) &= x + \frac{cx}{x - 1}, \\
    F_c(x) &= h(\psi_c(x)).
    \tag{$\star$}
\end{align*}

\subsubsection{Case: \texorpdfstring{$0<c\leq 1$}{0c1}}
Assume that $\lambda_k > \bbp$. Note that $h$ is strictly increasing on $[1, \infty)$, and $\psi_c$ is strictly increasing on $[\bbp, \infty)$, with $\psi_c(\bbp) = (\bbp)^2 =: b$. Therefore $F_c$ is also strictly increasing on $[\bbp, \infty)$.

\begin{figure}[!htbp]
    \centering
    \includegraphics[scale = 0.50]{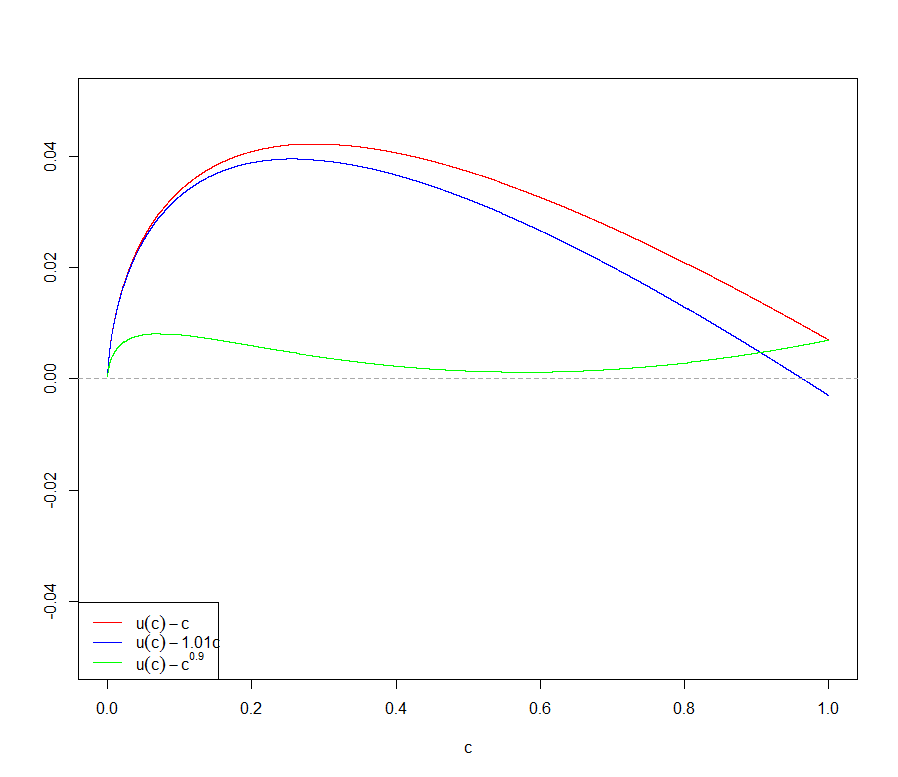}
    \caption{The size of the gap $u(c)$.}
    \label{fig:gapsize}
\end{figure}

Next we break down the proof of consistency by \cite{bai2018consistency} in two main steps. This gives us the important clue on how to modify the penalty of AIC.
\begin{itemize}
    \item[(i)] For showing $\hat{k}_{A} \le k$ a.s., \cite{bai2018consistency} needed the condition that
    \begin{align*}
        1 - b + \log b + 2 c > 0,    
    \end{align*}
    i.e. $2c > F_c(\bbp)$, which is true for any $c \in (0, 1)$ (as can be verified directly).
    \item[(ii)] To ensure that $\hat{k}_{A} \ge k$ a.s. \cite{bai2018consistency} needed the GAP condition (\ref{eq:C3}):
    \begin{align*}
    \psi_k - 1 - \log \psi_k - 2c > 0 \Leftrightarrow\,\, & F_c(\lambda_k) > 2c \\
    \Leftrightarrow\,\, &\lambda_k > F_c^{-1}(2c) = \psi_c^{-1}(h^{-1}(2c)) =: \lambda_c.  
    \end{align*}
\end{itemize}
Clearly $\lambda_c > \bbp$. Let 
\[
    u(c) := \lambda_c - (\bbp).  
\]
Let $u(c)$ denote the gap between the BBP threshold $\bbp$ and the consistency threshold $\lambda_c$ for AIC.
A natural question is, how large is the gap as a function of $c$? It is displayed in Figure~\ref{fig:gapsize}. The numerical calculation below shows that $u(c) > c$ ; in fact, $u(c) > c^{0.9}$.  

Our aim here is to investigate if we can develop a criterion that will be strongly consistent even when $\lambda_k$ is closer to $1+ \sqrt{c}$ than $\lambda_c$. The inequalities in (i) and (ii) above appear in the final asymptotic limits in the argument of \cite{bai2018consistency}. Placed together, for every $0<c \leq 1$, it is required that the inequalities $F_c(\lambda_k) > 2c > F_c(1+\sqrt{c})$ hold, where the factor 2 appears unchanged from the per parameter penalty 2 in AIC and $c$ from the limiting ratio of $p$ and $n$. Our modified AIC criterion is based on the crucial analytic observation we now make. From the monotonicity of of $F_c(.)$, the second inequality above also holds if 2 is replaced by any $\alpha < 2$ such that $F_c(1+\sqrt{c}) > \alpha c$. Evidently, for any fixed $c \in (0,1]$, there is whole interval of such choices of $\alpha$, and for that matter, one can choose $\alpha$ to make $F_c(1+\sqrt{c})$ arbitrarily close of $\alpha c$. For any such $\alpha$, a modified first inequality $F_c(\lambda_k) < \alpha c$ will be satisfied by any $\lambda_k$ larger than $F_c^{-1}(\alpha c))$, the latter being smaller than $\lambda_c$ by monotonicity of $F_c(.)$. Clearly $F_c^{-1}( \alpha c) - (1+\sqrt{c}) < u_c$. Smaller the choice of $\alpha$, closer is $F_c^{-1}( \alpha c)$ to $(1+\sqrt{c})$. We derive our modified AIC criterion by formally using this basic observation as described below.\\

Suppose we modify the original penalty of AIC with $\alpha(p/n)$, for some continuous function $\alpha:\mathbb{R} \mapsto (0, \infty)$, then the same arguments as in \cite{bai2018consistency} lead to the following conditions:
\begin{itemize}
    \item[(iii)] For showing $\hat{k}_{\alpha} \le k$ a.s.,
    \[
        1 - b + \log b + \alpha(c) c > 0,
    \]
    i.e. $\alpha(c) c > F_c(\bbp)$.
    \item[(iv)] For showing $\hat{k}_{\alpha} \ge k$ a.s.,
    \[            
        F_c(\lambda_k) > \alpha(c)c,
    \]
\end{itemize}
where $\hat{k}_{\alpha}$ is the minimizer of this new criterion over the model space.
 Thus if we want to get a modified criterion that is strongly consistent for any ``gap'' $\delta_c \in (0, u(c))$, all we need is that 
\[
    \alpha(c) = \frac{1}{c}F_c(\bbp + \delta_c).
\]
The monotonicity of $F_c(.)$ will ensure that the ineualities in (iii) and (iv) will be satisfied for this choice of $\alpha(c)$ which assures strong consistency of $\hat{k}_\alpha$ as long as $\lambda_k > \bbp + \delta_c$.
Thus we can reach arbitrarily close to the BBP threshold $\bbp$ and still achieve strong consistency. In fact, since $u(c) > c$, we may choose $\delta_c \in (0, c]$. For example, if we choose $\delta_c = c^{100} \ll c$, then with
    \[
        \alpha(c) = \frac{1}{c}F_c(\bbp + c^{100}),  
    \] 
we have consistency of $\hat{k}_{\alpha}$ as long as $\lambda_k > \bbp + c^{100}$. This is a better estimator than AIC because it can consistently estimate $k$ even when GAP condition (\ref{eq:C3}) fails to hold. In particular if we want our gap between $\lambda_k$ and $\bbp$ to be $\delta$, where $\delta>0$ is any arbitrarily small constant we can choose $\alpha(.)$ such that 
\[
        \alpha(c) = \frac{1}{c}F_c(\bbp + \delta),  
\]
and we would have strong consistency of our estimator as long as $\lambda_k > \bbp+\delta$.

We now describe in detail the new model selection criterion which we discussed above. The model selection criterion being a modification of AIC is being called $\aic^*$. Do describe this, we first fix a constant $\delta>0$.
We shall discuss the $p<n$ case first. With $C_{p,n}$ as defined before, 
the criterion value for $\aic^*$ for the model $M_j$ is defined as
\[
\aic^*_j = n \log (l_{1}\hdots l_{j}) + n(p-j) \log \bar{l}_{j} + \frac{1}{p/n}F_{p/n}(1+\sqrt{p/n} + \delta) d_j + C_{p,n},
\]
where, as before,  $l_{1}>\hdots >l_{p} > 0$ are the sample eigenvalues of $\mathbf{S_n}$ and for $1\leq j \leq p-1$, $\bar{l}_{j}$ is the arithmetic mean of $l_{j+1},\hdots, l_{p}$, that is,
\[
\bar{l}_{j} = \frac{1}{p-j}\sum_{t = j+1}^p l_{t}  
\]
and, $d_j =(j+1)(p+1-j/2)$ denotes the number of independent parameters in the model under $M_j$.
Furthermore the function $F_c(.)$ is defined as in (\ref{eq:*}).
Then the $\aic^*$ selects $k$ according to the rule
\[
\hat{k}^*_A = \text{arg min}_j \aic^*_j,
\]
where the minimum has been taken over the model space.

Now we state our main result regarding consistency of our estimator $\hat{k}^*_A$ whose proof can be found in Appendix Section \ref{strongconsproof1}. 
\begin{theorem}
\label{strongconsistencythm1}
Suppose the conditions (\ref{eq:C1}) with $0<c \leq 1$, and (\ref{eq:C2}) hold, and that the number of candidate models, $q$, satisfies $q=o(p)$. Assume that $\lambda_1$ is bounded. We have the following results on the consistency of the estimation criterion $\hat{k}^*_A$
based on $\aic^*$:\\

\begin{itemize}
    \item[(i)] If $\lambda_k \leq \bbp +\delta$, then $\hat{k}^*_A$ is not consistent.
    \item[(ii)] If $\lambda_k > \bbp +\delta$, then $\hat{k}^*_A$ is strongly consistent.
\end{itemize}
\end{theorem}

\subsubsection{Case: \texorpdfstring{$c>1$}{c1}}
Next we consider the case where $c>1$.
Let
\[
\label{eq:**}
Q_c(x) = c h\left(\frac{\psi_c(x)}{c}\right).
\tag{$\star\star$}
\]
 Recall that $h$ is strictly increasing on $[1, \infty)$, and $\psi_c$ is strictly increasing on $[\bbp, \infty)$, {\it Therefore} $\psi_c(x)/c$ is increasing in $x$. We have $\psi_c(\bbp)/c = (\bbp)^2/c > 1 \quad \forall c>1$. Therefore $Q_c$ is also strictly increasing on $[\bbp, \infty)$.
\begin{itemize}
    \item[(i)] For showing $\hat{k}_{\tilde{A}} \le k$ a.s., \cite{bai2018consistency} needed the condition
    \begin{align*}
        1 - b/c + \log (b/c) + 2/c > 0,    
    \end{align*}
    i.e. $2>Q_c(\bbp)$, which is true for any $c >1$ (as can be verified directly).
    \item[(ii)] To ensure that $\hat{k}_{\tilde{A}} \ge k$ a.s. \cite{bai2018consistency} need the GAP condition (\ref{eq:C5}):
    \begin{align*}
    \psi_k/c - 1 - \log (\psi_k/c) - 2/c > 0 \Leftrightarrow
     Q_c(\lambda_k) > 2
    \Leftrightarrow \lambda_k > Q_c^{-1}(2) =: \lambda_c.  
    \end{align*}
\end{itemize}
 
\begin{figure}[!htbp]
    \centering
    \includegraphics[scale = 0.60]{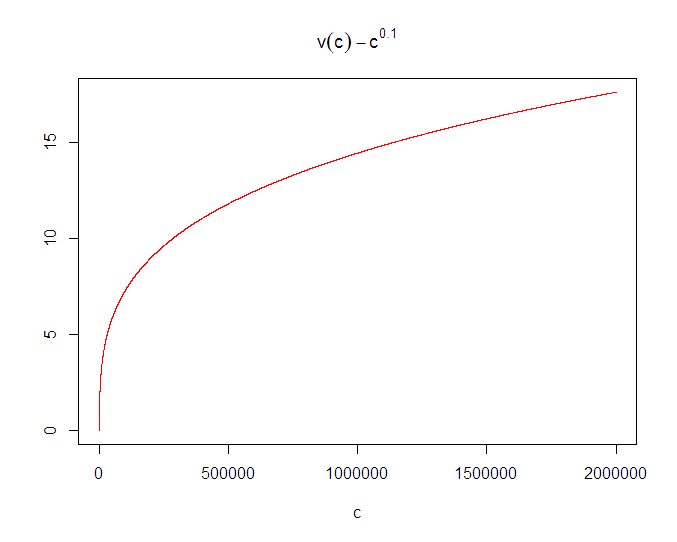}
    \caption{The size of the gap $v(c)$.}
    \label{fig:gapsize_1}
\end{figure}

Note that $\lambda_c > \bbp$. Let 
\[
    v(c) := \lambda_c - (\bbp).  
\]
$v(c)$ is the gap between the {\it BBP} threshold $\bbp$ and the consistency threshold $\lambda_c$ for qAIC.
Our aim  as in the previous case is to {\it investigate if a new criterion can be derived with is strongly consistent even when the ``gap'' is reduced} as much as possible. So firstly we look at how large is the gap as a function of $c$? Figure~\ref{fig:gapsize_1} shows that $v(c) > c^{0.1}$.  

Following similar intuition as in the case of $0<c \leq 1$, if we penalize by replacing the $2$ in qAIC with $\alpha(p/n)$ for some function $\alpha:\mathbb{R} \mapsto (0, \infty)$, then the same arguments as in \cite{bai2018consistency} lead to the following conditions:
\begin{itemize}
    \item[(iii)] For showing $\hat{k}_{\alpha} \le k$ a.s.,
    \[
        1 - b/c + \log (b/c) + \alpha(c)/c > 0,
    \]
    i.e. $\alpha(c) > Q_c(\bbp)$.
    \item[(iv)] For showing $\hat{k}_{\alpha} \ge k$ a.s.,
    \[            
        Q_c(\lambda_k) > \alpha(c).
    \]
\end{itemize}

Thus for any $\delta_c \in (0, v(c))$, if we take the function 
\[
    \alpha(c) =Q_c(\bbp + \delta_c),
\]

then using the monotonicity of $Q_c(.)$ we have consistency of $\hat{k}_\alpha$ as long as $\lambda_k > \bbp + \delta_c$. Thus we can reach arbitrarily close to the {\it BBP} threshold $\bbp$.

Our model selection criterion is {\it a} modification of the criterion proposed by \cite{bai2018consistency} as discussed in Section \ref{sec:baistrongconsistency}. We define our criterion using $\tilde{A^*}$, which for model $M_j$ is given by
\[
\tilde{A^*_j} = (n-1-j)\log \bar{l}_{j} - \sum_{i=j+1}^{n-1}\log l_{i} - Q_{p/n}(1+\sqrt{p/n}+\delta)\frac{(n-j-2)(n-j+1)}{2p}
\]
where $\bar{l}_{jp}$ is defined as in \cite{bai2018consistency} for $p>n$ case which is
$$\bar{l}_{j} = \frac{1}{n-i-1}\sum_{t= i+1}^{n-1} l_{t}$$
and $Q_c(.)$ is as defined in (\ref{eq:**}).
The model selection rule $\tilde{A^*}$ estimates $k$ by $ \hat{k}^*_{\tilde{A}}$ defined as
\[
\hat{k}^*_{\tilde{A}} = \text{arg min}_j \tilde{A^*}_j,
\]
 the minimum. as before, being taken over the model space.

The following result shows the consistency of our estimator $\hat{k}^*_{\tilde{A}}$. 
\begin{theorem}
\label{strongconsistencythm2}
Suppose the conditions (\ref{eq:C1}) with $c>1$ and (\ref{eq:C4}) hold, and that the number of candidate models, $q$, satisfies $q=o(p)$. Suppose that $\lambda_1$ is bounded. We have the following results on the consistency of the estimation criterion $\hat{k}^*_{\tilde{A}}$:
\begin{itemize}
    \item[(i)] If $\lambda_k \leq \bbp +\delta$, then $\hat{k}_{\tilde{A}}$ is not consistent.
    \item[(ii)] If $\lambda_k > \bbp +\delta$, then $\hat{k}_{\tilde{A}}$ is strongly consistent.
\end{itemize}
\end{theorem}

\subsection{Weak Consistency for a modified AIC criterion}
\label{sec:weakconsistencyresult}
We have seen in the previous section that for any given $\delta > 0$, however small, we can appropriately modify the penalty term (using $\delta$ among other things) of the usual AIC to produce an estimator of $k$ which is strongly consistent when the gap between $\lambda_k$ and the BBP threshold is larger than $\delta$. But ideally we would aim to develop consistent estimators of $k$ which would be consistent if $\lambda_k > 1+\sqrt{c}$, i.e. the we can afford to have the ``gap'' equal to zero. A natural thing is to first investigate what happens when we choose the gap to be $\delta_n>0$, where $\delta_n \to 0$ as $n \rightarrow \infty$ in the definition of the modified AIC criterion described before.

We first discuss the case $0 <c \leq 1$. We mofify the penalty term in AIC criterion as $\alpha(p,n)$ instead of $2$ where
$$\alpha(p,n) = \frac{1}{p/n}F_{p/n}(1+\sqrt{p/n}+\delta_n).$$ We define $\hat{k}_{\alpha}$ as the index that minimizes this modified criterion. Observe that $$\alpha(p,n) \to \frac{1}{c}F_c(\bbp)\qquad \text{ as } n\to \infty.$$ Let us check if $\hat{k}_{\alpha}$ is strongly consistent. Following the arguments of \cite{bai2018consistency}, we find the following.  
\begin{itemize}
    \item[(i)] For showing $\hat{k}_{\alpha} \ge k$ a.s. we need,
    \[            
        F_c(\lambda_k) > F_c(\bbp).
    \]
    This is true since we have that $\lambda_k > \bbp$.
    \item[(i)] For showing $\hat{k}_{\alpha} \le k$ a.s. we need,
    \[
        1 - b + \log b + F_c(\bbp) > 0,
    \]
    i.e. $F_c(\bbp) > F_c(\bbp)$, which is not true.
\end{itemize}

Clearly the arguments given by \cite{bai2018consistency} doesn't give us our desired result for the above intuitive modification of our previous criterion. We however have been able to come up with a novel argument that gives us our desired result. The argument is described in detail in Section \ref{weakconsproof1} of the Appendix and carefully exploits some deep results available from Random Matrix Theory regarding the rate of convergence of certain functionals of the sample covariance matrix. Ultimately it turns out that $\delta_n$ has to converge to zero in such a rate that the new penalty $\alpha(p,n)$ satisfies certain convergence properties. The only downside is that we cannot prove strong consistency of the new criterion using this method, but we can show weak consistency of the estimator under some assumptions (see Theorem \ref{weakconsitencythm1} and \ref{weakconsistencythm2}).

As mentioned above, the basic idea for developing the estimator which would be consistent under the ``zero'' gap condition was to use our previous modified AIC criterion with $\delta_n$ instead of $\delta$, where $\delta_n \to 0$. We describe below the model selection rule. This The being a modification of $\aic^*$, we call it $\aic^{**}$.
We discuss the $p<n$ case first. With  $C_{p,n}$ defined as before, the criterion value under model $M_j$ is defined as
\[
\aic^{**}_j = n \log (l_{1}\hdots l_{j}) + n(p-j) \log \bar{l}_{j} + \frac{1}{p/n}F_{p/n}(1+\sqrt{p/n} + \delta_n) d_j + C_{p,n}
\]
Furthermore, $d_j$ is as defined before.
Defining
$$A^{**}_j = \frac{1}{n}(\aic^{**}_j -\aic^{**}_{p-1}),$$
we have
\[
A^{**}_j = (p-j)\log \bar{l}_{j} - \sum_{i=j+1}^p \log l_{i} -\frac{1}{p/n}F_{p/n}(1+\sqrt{p/n} + \delta_n) \frac{(p-j-1)(p-j+2)}{2n}
\]
Then the $\aic^{**}$ criterion selects model with index  $\hat{k}^{**}_A$ obtained as
\[
\hat{k}^{**}_A = \text{arg min}_j \aic^{**}_j,
\]
or equivalently
\[
\hat{k}^{**}_A = \text{arg min}_j A^{**}_j.
\]
The criteria above is defined by considering the minimum only with respect to $j=0,1,\hdots,s-1$.
where $s > k$ is a fixed number big enough. In practice, the integer
$s$ should be thought as a preliminary bound on the number of
possible spikes.
We have the following theorem regarding consistency of the proposed estimator whose proof can be found in Appendix Section \ref{weakconsproof1}.
\begin{theorem}
\label{weakconsitencythm1}
Suppose the condition (\ref{eq:C1}) holds for $0<c\leq 1$ and (\ref{eq:C2}) hold, and that that the number of candidate models is $s$, where $s>k$ is a fixed number. Suppose that \hyperlink{assump1}{\text{Assumption 1}} holds true, $\lambda_1$ is bounded and $\lambda_k > \bbp$. Let $\delta_n$ be a real sequence going to 0 such that 
$$n^{2/3}\left(F_{p/n}(1+\sqrt{p/n}+\delta_n)  -F_c(1+\sqrt{c})\right) \to \infty,$$
where $F_c$ is as defined in (\ref{eq:*}),
then the estimator $\hat{k}^{**}_A$ is a weakly consistent estimator of $k$, i.e. $\mathbb{P}(\hat{k}^{**}_A = k)\to 1$ as $n \to \infty$. 
\end{theorem}

We have the similar result for the case $p>n$, i.e. $c>1$. Our modified model selection criterion is defined using $\tilde{A_j}^{**}$, which is defined as
\[
\tilde{A_j}^{**} = (n-1-j)\log \bar{l}_{j} - \sum_{i=j+1}^{n-1}\log l_{i} - Q_{p/n}(1+\sqrt{p/n}+\delta_n)\frac{(n-j-2)(n-j+1)}{2p}
\]
Then $\tilde{A_j}^{**}$  selects the number of significant components according to
\[
\hat{k}^{**}_{\tilde{A}} = \text{arg min}_j \tilde{A_j}^{**}. 
\]
The criteria above is defined by considering the minimum only with respect to $j=0,1,\hdots,s-1$.
where $s > k$ is a fixed number big enough. In practice, the integer
$s$ should be thought as a preliminary bound on the number of
possible spikes.
We have the following theorem regarding consistency of the proposed estimator whose proof can be found in Appendix Section \ref{weakconsproof2}.
\begin{theorem}
\label{weakconsistencythm2}
Suppose the condition (\ref{eq:C1}) with $c>1$ and (\ref{eq:C2}) hold, and that that the number of candidate models is $s$, where $s>k$ is a fixed number. Suppose that \hyperlink{assump1}{\text{Assumption 1}} holds true, $\lambda_1$ is bounded and $\lambda_k > \bbp$. We also assume that $p/n = c + O(n^{-2/3})$.\\
Let $\delta_n$ be a real sequence going to 0 such that 
$$n^{2/3}\left( \frac{1}{p/n}Q_{p/n}(1+\sqrt{p/n}+\delta_n) -\frac{1}{c}Q_c(\bbp)\right) \to \infty,$$
where $Q_c$ is as defined in (\ref{eq:**}),
then the estimator $\hat{k}^{**}_{\tilde{A}}$ is a weakly consistent estimator of $k$, i.e. $\mathbb{P}(\hat{k}^{**}_{\tilde{A}} = k)\to 1$ as $n \to \infty$. 
\end{theorem}

\section{Simulation studies}
\label{sec:simulation}

\subsection{The choice of \texorpdfstring{$\boldsymbol{\delta_n}$}{d}}
\label{sec:choice of delta_n}
Consistency of our estimator holds for a wide variety of sequences. For instance, if $\delta_n \to 0$ is a sequence satisfying the condition of our theorem, then so does any $c \delta_n$ where $c \in \mathbb{R}$ and $c>0$ . In this section we provide an automatic calibration of this parameter. The main idea is to look at the null case ($k=0$), i.e. $X_i \overset{i.i.d.}{\sim} N(0,I_p)$. Given a $n$ and a $p$, we propose to select $\delta_n$ as 
$$\delta_n = \inf\{\delta>0 \mid \mathrm{RMSE}(\hat{k}_{\delta,n}) \leq 0.02\},$$
where $\hat{k}_{\delta,n}$ is our strongly consistent estimator dependent on $\delta$ as proposed in Theorem \ref{strongconsistencythm1} or \ref{strongconsistencythm2} and RMSE or Relative Mean Square Error is given by
$$\mathrm{RMSE}(\hat{k}_{\delta,n}) = \mathbb{E}\bigg(\frac{\hat{k}_{\delta,n}-k}{k}\bigg)^2.$$
As we do not know the precise expression of RMSE for the null case we approximately calculate it by simulation under $500$ independent replications and call it SRMSE (Sample RMSE). Using the SMRSE in place of RMSE above we find our choice of $\delta_n$.   
 
Table~\ref{table:choice_of_delta_n} lists the values of $\delta_n$ chosen according to our proposal for a few different values of $(p,n)$.
\begin{table}[!h]
\centering
\begin{tabular}{|c|c|c|c|c|c|c|}
\hline
$(p,n)$    & (200,200) & (400,400) & (1000,1000) & (200,400) & (500,1000) & (400,500) \\ \hline
$\delta_n$ & 0.52      & 0.34      & 0.25        & 0.28      & 0.20       & 0.32      \\ \hline
\end{tabular}
\caption{The choice of $\delta_n$ for different values of $(p, n)$.}
\label{table:choice_of_delta_n}
\end{table}
Let us now illustrate with the help of two models how our method of choosing $\delta_n$ compares with manually choosing $\delta_n$ (i.e. by choosing the $\delta_n$ which gives the ``lowest'' RMSE).

$X_1,X_2,\hdots,X_n$ are generated from $N_p(0,\mathbf{\Sigma})$ where $\mathbf{\Sigma} = \text{diag}\{\lambda_1,\hdots,\lambda_k,1,\hdots,1\}$.
 After fixing a $c$, $n$ and $p$ are varied such that $p/n \approx c$.
 We compute the RMSE (by $400$ replications) for our estimator, for the two different choices of $\delta_n$, first one is automatic calibration, and the second one by manual tuning (to get the best performance). We compare both these method for the following two models.
 
\begin{figure}[!h]
    \centering
    \begin{tabular}{cc}
    \qquad$c=1,\, k=2,\, \lambda_1=4.5, \lambda_2 = 3$ & \qquad$c=1,\, k=2,\, \lambda_1 = 3, \lambda_2 =2.3$\\
    \includegraphics[scale = 0.55]{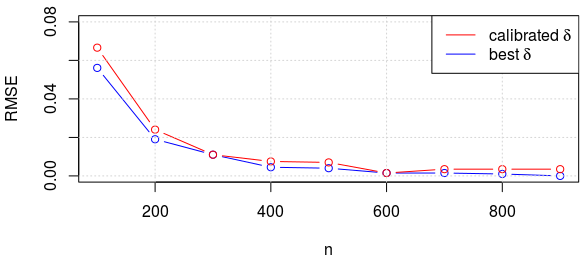}
    & 
    \includegraphics[scale =0.55]{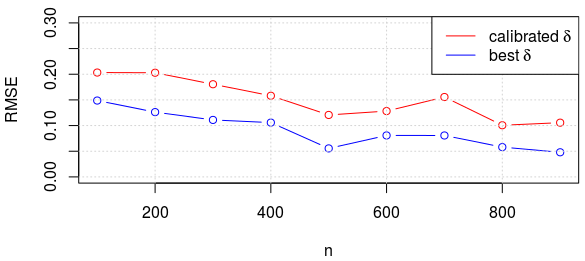} \\
    \qquad\,\,\,(a)
    & 
    \qquad\,\,\,(b)
    \end{tabular}
    \caption{(a) Model 1; (b) Model 2.}
    \label{fig:choice_of_delta_n}
\end{figure}

\begin{itemize}
    \item \textbf{Model 1:} $c=1,\, k=2,\, \lambda_1=4.5, \lambda_2 = 3$. In this case our calibrated $\delta_n$ has performance really close to the ``ideal'' $\delta_n$ as can be seen by the Figure \ref{fig:choice_of_delta_n}-(a).

    \item \textbf{Model 2:} $c=1,\, k=2,\, \lambda_1 = 3, \lambda_2 =2.3$. In this case the gap between BBP threshold and $\lambda_k$ is very small $(=0.3)$. Here the difference between the two choice of $\delta_n$ is more prominent. The RMSE corresponding to the ``best'' $\delta_n$ is about $25\%$ lower than that of our calibrated $\delta_n$ (See Figure \ref{fig:choice_of_delta_n}-(b)).
\end{itemize}

The main takeaway here is that when the gap is not too small then our calibrated $\delta_n$ has performance very close to the ``best'' $\delta_n$, but when the ``gap'' is very small, further study needs to be done to be able to adaptively choose $\delta_n$ based on the ``gap''. Nevertheless in either case (gap: small/big) our modified AIC criterion (with calibrated $\delta_n$) works better than the other estimators in the literature as shown in the next two sections.
\subsection{Comparison study I}
In this section we compare our modified AIC estimator against the estimator by  \cite{passemier2014estimation} (abbreviated as PY), see section \ref{sec:Passemier and Yao} for details.

\begin{figure}[!h]
    \centering
    \begin{tabular}{cc}
        \qquad$c =1$, $k=2$, $\lambda_1 = 3.5, \lambda_2 =2.5$ & \qquad$c=1$, $k=2$, $\lambda_1 = 3, \lambda_2 =2.1$ \\
        \includegraphics[scale =0.55]{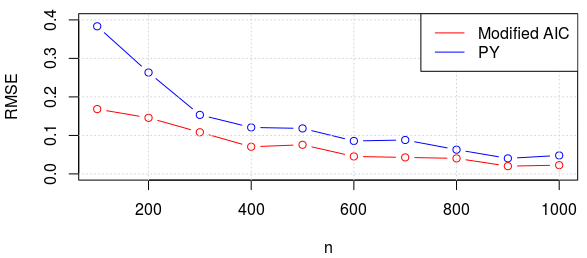}
        &
        \includegraphics[scale =0.55]{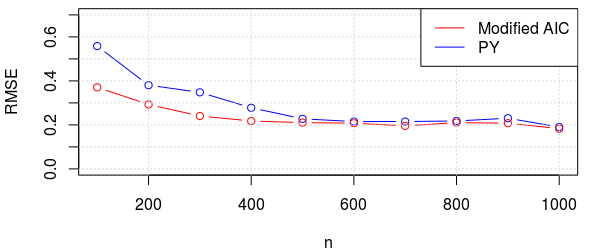} \\
        \qquad\,\,\,(a) & \qquad\,\,\,(b) \\
         & \hskip-250pt \qquad$c =1$, $k=2$, $\lambda_1 = 3, \lambda_2 =3$ \\
         & \hskip-250pt\includegraphics[scale =0.55]{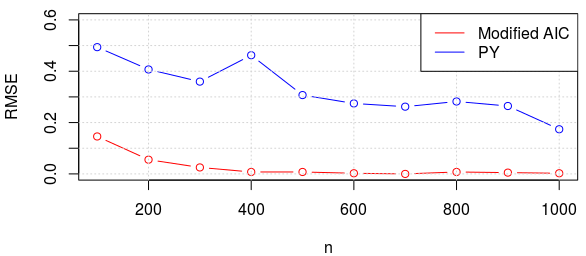} \\
         & \hskip-250pt\qquad\,\,\,(c)
    \end{tabular}
    \caption{(a) Model A; (b) Model B; (c) Model (c).}
    \label{fig:comp_study_1}
\end{figure}

$X_1,X_2,\hdots,X_n$ are generated from a $p$-variate Gaussian distribution with mean 0 and a spiked variance-covariance matrix $\Sigma = \text{diag}\{\lambda_1,\hdots,\lambda_k,1,\hdots,1\}$. 
 After fixing a $c$, $n$ and $p$ are varied such that $p/n \approx c$.
 We compute the RMSE for our estimator (with calibrated $\delta_n$) and for PY estimator (as proposed in their paper with automatic calibration) by conducting $100$ independent replications. We compare both these method for the following models.
\begin{itemize}
    \item \textbf{Model A:} $c =1$, $k=2$, $\lambda_1 = 3.5, \lambda_2 =2.5$.
    Our method performs uniformly better than PY (see Figure~\ref{fig:comp_study_1}-(a)), though the advantage seems to decrease as $n$ increases.

    \item \textbf{Model B:} $c=1$, $k=2$, $\lambda_1 = 3, \lambda_2 =2.1$. In this case the gap between BBP threshold and $\lambda_k$ is very small $(=0.1)$. Our method is still the better one (see Figure~\ref{fig:comp_study_1}-(b)).

    \item \textbf{Model C:} $c =1$, $k=2$, $\lambda_1 = 3, \lambda_2 =3$ (equal spikes).
    Our Method is clearly better than PY's (see Figure~\ref{fig:comp_study_1}-(c)) and the advantage gained is significant especially in case of equal spikes.
\end{itemize}

\subsection{Comparison study II}
In this section we compare our estimator against the ordinary AIC criterion as proposed by \cite{bai2018consistency}, see section \ref{sec:baistrongconsistency} for details. 

\begin{figure}[!h]
    \centering
    \begin{tabular}{cc}
    \qquad$c =1$, $k=2$ $\lambda_1 = 3.5, \lambda_2 =3.5$ & \,\,\,\,$c =1$, $k=5$ $\lambda_1 = 4, \lambda_2 =4, \lambda_3=\lambda_4=\lambda_5 =3.5$ \\
    \includegraphics[scale =0.55]{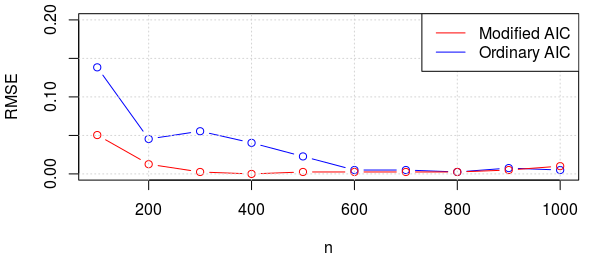}
    &
    \includegraphics[scale =0.55]{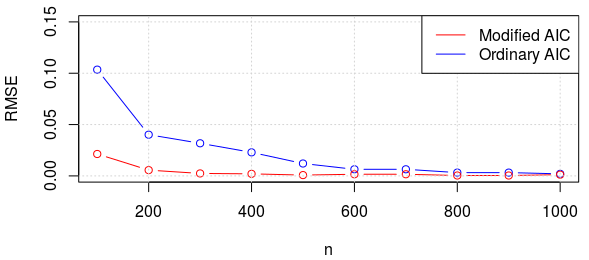} \\
    \qquad\,\,\,(a) & \qquad\,\,\,(b)
    \end{tabular}
    \caption{(a) Model D; (b) Model E.}
    \label{fig:comp_study_2}
\end{figure}

$X_1,X_2,\hdots,X_n$ are generated from a $p$-variate Gaussian distribution with mean 0, with spiked variance-covariance matrix $\Sigma = \text{diag}\{\lambda_1,\hdots,\lambda_k,1,\hdots,1\}$. 
 After fixing a $c$, $n$ and $p$ are varied such that $p/n \approx c$.
 We compute the RMSE for our estimator (with calibrated $\delta_n$) and for ordinary AIC criterion by conducting $100$ independent replications. We compare both these method for the following models.
\begin{itemize}
    \item \textbf{Model D:} $c =1$, $k=2$ $\lambda_1 = 3.5, \lambda_2 =3.5$.
    Here $\lambda_k =\lambda_2$ is chosen such that it satifies the gap condition (\ref{eq:C3}) as proposed by  \cite{bai2018consistency}.
    Observe that our method performs uniformly better compared to theirs (see Figure~\ref{fig:comp_study_2}-(a)), though the advantage seems to decrease as $n$ increases.

    \item \textbf{Model E:} $c =1$, $k=5$ $\lambda_1 = 4, \lambda_2 =4, \lambda_3=\lambda_4=\lambda_5 =3.5$.
    Here $\lambda_k =\lambda_5$ is chosen such that it satisfies the gap condition (\ref{eq:C3}).
    Observe that in this case too our method performs uniformly better compared to theirs (see Figure~\ref{fig:comp_study_2}-(b)).
\end{itemize}

\subsection{The choice of \texorpdfstring{$\boldsymbol{\delta = 0}$}{d}}
\label{sec:delta as zero0}
Observe that we can reach arbitrary close to $\bbp$ , by suitably choosing the penalty term in AIC criterion ($\alpha(c)$).
Motivated by our estimators where we have chosen our gap to be $\delta$ where $\delta$ can be arbitrarily small, or even $\delta_n \to 0$ (in case of weak consistency), we present an estimator where $\delta$ is exactly zero, i.e. choose
\[
    \alpha(c) =\frac{1}{c} F_c(\bbp).
\]
We now present a simulation study about this estimator.
\begin{itemize}
\item \textbf{Simulation setup:}
$c = 0.5$ is fixed, $n$ is varied, $X_1,X_2,\hdots,X_n$ are generated from a $p$-variate Gaussian distribution with mean $0$ with spiked variance-covariance matrix (i.e. a diagonal matrix with first $k$ entries as $\bbp +\delta$ and rest $p-k$ diagonal entries are 1). Here we fixed $\delta =0.5$ and $k =10$. 
\item \textbf{Results:} We varied $n$ from $100$ to $3200$ and iterated $50$ times each. Table~\ref{table:delta_0} shows the accuracy (i.e. proportion of times $\hat{k}$ is equal to true $k$) and the average value of $\hat{k}$ over $50$ iterations.
\begin{table}[!h]
\centering
\begin{tabular}{|c|c|c|c|c|c|c|c|c|c|}
\hline
$n$             & 100  & 200  & 400  & 800  & 1200 & 1500 & 2000  & 2500  & 3200  \\ \hline
accuracy      & 0.02 & 0.04 & 0.04 & 0.26 & 0.48 & 0.76 & 0.94  & 0.96  & 0.94  \\ \hline
avg $\hat{k}$ & 5.56 & 7.08 & 7.82 & 8.98 & 9.50 & 9.82 & 10.06 & 10.00 & 10.06 \\ \hline
\end{tabular}
\caption{The behaviour of our estimator with the choice $\delta = 0$.}
\label{table:delta_0}
\end{table}
\end{itemize}
The results suggest that the estimator is consistent. We leave the theoretical analysis of this case to future work.

\section{Conclusion}
In this article, we have considered the AIC cirterion for high dimensional model selection. \cite{bai2018consistency} have shown strong consistency of AIC under a ``gap'' condition, requiring more signal than what the BBP threshold demands. We have modified the penalty term in AIC suitably and proposed a criterion that achieves strong consistency under arbitrarily small gap. We have also proposed another modification that achieves weak consistency under exactly zero gap. We have made a detailed empirical study of the performance of the proposed criteria. Furthermore, we have compared our second proposal with the estimator of \cite{passemier2014estimation} which also achieves weak consistency under zero gap.

\section{Acknowledgements}
This article is based on the Masters dissertation \cite{chakraborty2020thesis} of the first author under the supervision of the second and the third authors. The authors thank Professors Gopal K. Basak and Tapas Samanta for their insightful comments. The second author is supported by an INSPIRE Faculty Fellowship, Department of Science and Technology, Government of India.

\bibliographystyle{apalike}  
\bibliography{references}

\newpage
\section{Appendix}
\subsection{Preparatory Lemmas}
The proof of Theorem \ref{strongconsistencythm1}, \ref{weakconsitencythm1} and \ref{weakconsistencythm2} rely on several additional lemmas which are organised in this section. We first recall a lemma stated earlier in Section~\ref{sec:rmt_results} as Lemma~\ref{asconvergence}.
\begin{lemma}
Let $l_{ip}$ denote the ith largest eigenvalue of $S_n$ the covariance matrix in our setup. Suppose that $\mathbb{E}(x_{11}^4 <\infty )$, condition $(\ref{eq:C1})$ and $(\ref{eq:C2})$ hold, and that $\lambda_1$ is bounded.
\begin{itemize}
    \item[(i)] If $\lambda_i > 1+\sqrt{c}$, then $l_{ip} \overset{a.s.}{\to}\psi_i = \lambda_i + \frac{c\lambda_i}{\lambda_i -1}$.
    \item[(ii)] $\lambda_i < 1+\sqrt{c}$ and $i/p \to \alpha$, then $l_{ip}\overset{a.s.}{\to} \mu_{1-\alpha}$(where $\mu_{\alpha}$ is the $\alpha$-quantile of MP distribution). In particular if i is bounded $l_{ip}\overset{a.s.}{\to} \mu_1 =b:= (1+\sqrt{c})^2$.  
\end{itemize}
\end{lemma}

The next lemma is a direct consequence of Theorem 3.1 of, \cite{li2019asymptotic}.
\begin{lemma}
\label{LiHanYao}
Under the given assumptions of the Theorem \ref{weakconsitencythm1} or \ref{weakconsistencythm2} we have that 
$$trace(S_n)-trace(V) \overset{d}{\to} N(0,\sigma_0^2),$$
where $\sigma_0$ is some constant independent of $n$ or $p$.
\end{lemma}

Next lemma is a direct consequence of Proposition $5.8$ from  \cite{benaych2011fluctuations}.
\begin{lemma}
\label{passemierlemma}
Assume that the entries $x^i$ of $x$ have a symmetric law and a sub-exponential decay, that means $\exists$ positive constants $C, C'$ such that, for all $t>C'$, we have $\mathbb{P}(|x^i| \geq t^C) \leq e^{-t}$. Then, for all $1\leq i \leq L$ with a prefixed range L,
$$
n^{\frac{2}{3}}(\lambda_{k+i} -b) = O_{\mathbb{P}}(1).
$$
\end{lemma}

\subsection{Weak Consistency Proof}

In this section we are going to prove Theorem \ref{weakconsitencythm1} and \ref{weakconsistencythm2}.

\label{weakconsproof1}
\textbf{Proof of the Theorem \ref{weakconsitencythm1}.} We will show that $\mathbb{P}(\hat{k}^{**}_A = k)\to 1$ as $n \to \infty$. We will show that by breaking the problem into two cases $j<k$ and $j>k$, in both these cases we will show that $\mathbb{P}(A^{**}_j >A^{**}_k)\to 1$ as $n \to \infty$.\\ 
 \textbf{Case I:} When $j<k$,
 \begin{align*}
     A^{**}_j-A^{**}_k &= \sum_{i=j+1}^k (A^{**}_{i-1} - A^{**}_i) \\
     \begin{split}
     &= \sum_{i=j+1}^k 
     (p-i+1)\log\left\{1-\frac{1}{p-i+1} \left(1-\frac{l_{ip}}{\bar{l}_{ip}}\right)  \right\}\\
     & +\log\bar{l}_{ip} - \log l_{ip} - \frac{1}{p/n}F_{p/n}(1+\sqrt{p/n}+\delta_n)\frac{p-i+1}{n} 
    \end{split}\\
    &\overset{a.s.}{\to} \sum_{i=j+1}^k \psi_i - 1 - \log\psi_i - F_c(\bbp) \\
    &= \sum_{i=j+1}^k F_c(\lambda_i) - F_c(\bbp) 
 \end{align*}
 Only step that needs to be justified is the third step , specifically\\
 \textbf{Claim:}
 $$(p-i+1)\log\left\{1-\frac{1}{p-i+1} \left(1-\frac{l_{ip}}{\bar{l}_{ip}}\right)\right\} \overset{a.s.}{\to} \psi_i -1.$$
Observe that,
     $\left|(p-i+1)\log\left\{1-\frac{1}{p-i+1} \left(1-\frac{l_{ip}}{\bar{l}_{ip}}\right)  \right\} -(\psi_i-1)\right|$ \\
     $\leq$ $\left|-(p-i+1)\log\left\{1-\frac{1}{p-i+1} \left(1-\frac{l_{ip}}{\bar{l}_{ip}}\right)  \right\} - \left(1-\frac{l_{ip}}{\bar{l}_{ip}}\right)\right|$ +$\left| \frac{l_{ip}}{\bar{l}_{ip}} - \psi_i \right|.$\\
     
     \textbf{Result:(Taylor series theorem)} For a twice differentiable $f$ on $(a,b)$ with $f'$ continuous on $[a,b]$. Assume that $c \in [a,b]$. Then for every $x \in [a,b] \neq c$, there exists point $x_1$ in $(\min(x,c),\max(x.c))$ such that
     $$ f(x) = f(c) + f'(c)(x-c) + \frac{f''(x_1)}{2}(x-c)^2.$$\\
     Using the above theorem we get under prescribed conditions on $f$ we have
     $$f(x+h_0) = f(x) + f'(x)h_0 + \frac{f''(\zeta)}{2}h_0^2 \quad where \,\, \zeta \in (\min(x,x+h_0),\max(x.x+h_0)).$$
     Choose 
     $$f(x) =\log(1-x)\qquad x=0  \qquad h_0 = \frac{1}{p-i+1}\left(1-\frac{l_{ip}}{\bar{l}_{ip}}\right)$$
     and observe that
     $$f'(x) = -\frac{1}{1-x} \qquad f''(x) = \frac{1}{(1-x)^2}.$$
     We have that
     $$\log\left\{1-\frac{1}{p-i+1} \left(1-\frac{l_{ip}}{\bar{l}_{ip}}\right)  \right\} = -\frac{1}{p-i+1}\left(1-\frac{l_{ip}}{\bar{l}_{ip}}\right)+ \frac{h_0^2}{2(1-\zeta)^2}$$
     where
     $\zeta \in (\min(0,h_0),\max(0,h_0))$.\\
     So we have
     $$-(p-i+1)\log\left\{1-\frac{1}{p-i+1} \left(1-\frac{l_{ip}}{\bar{l}_{ip}}\right)  \right\} = \left(1-\frac{l_{ip}}{\bar{l}_{ip}}\right)-(p-i+1) \frac{h_0^2}{2(1-\zeta)^2}.$$
     Hence
     \begin{align*}
     \left|-(p-i+1)\log\left\{1-\frac{1}{p-i+1} \left(1-\frac{l_{ip}}{\bar{l}_{ip}}\right)  \right\} - \left(1-\frac{l_{ip}}{\bar{l}_{ip}}\right)\right| &= \left|(p-i+1)\frac{h_0^2}{2(1-\zeta)^2}\right|\\
     =& \frac{1}{(1-\zeta)^2}\left(1-\frac{l_{ip}}{\bar{l}_{ip}}\right)^2\frac{1}{p-i+1}
     \end{align*}
     Lets look at the R.H.S, we know from Lemma \ref{asconvergence} that $l_{ip}\overset{a.s.}{\to}\psi_i$ and from the MP-law that $\bar{l}_{ip}\overset{a.s.}{\to} \mu_{MP} =1$, where $\mu_{MP}$ is the mean of the MP distribution.  \\
     Using these we infer that
     $$\left(1-\frac{l_{ip}}{\bar{l}_{ip}}\right)^2 \overset{a.s}{\to} (1-\psi_i)^2 \quad \& \quad h_0 = \frac{1}{p-i+1}\left(1-\frac{l_{ip}}{\bar{l}_{ip}}\right)\overset{a.s}{\to} 0. $$
     Now using sandwich theorem we have that $\zeta \overset{a.s}{\to} 0 $ as $\zeta \in (\min(0,h_0),\max(0,h_0))$ 
     so, we have that $\frac{1}{(1-\zeta)^2} \overset{a.s}{\to} 1$.\\
     Therefore,
     $$\frac{1}{(1-\zeta)^2}\left(1-\frac{l_{ip}}{\bar{l}_{ip}}\right)^2 \overset{a.s}{\to} (1-\psi_i)^2.$$
     So,
     $$\left|-(p-i+1)\log\left\{1-\frac{1}{p-i+1} \left(1-\frac{l_{ip}}{\bar{l}_{ip}}\right)  \right\} - \left(1-\frac{l_{ip}}{\bar{l}_{ip}}\right)\right| = \frac{\left(1-\frac{l_{ip}}{\bar{l}_{ip}}\right)^2}{(1-\zeta)^2}\frac{1}{p-i+1}$$
     As $\frac{1}{p-i+1} \to 0$ coupled with previous fact, we can conclude that 
     $$\left|-(p-i+1)\log\left\{1-\frac{1}{p-i+1} \left(1-\frac{l_{ip}}{\bar{l}_{ip}}\right)  \right\} - \left(1-\frac{l_{ip}}{\bar{l}_{ip}}\right)\right| \overset{a.s.}{\to} 0.$$
     Next using Lemma \ref{asconvergence} observe that $l_{ip} \overset{a.s.}{\to} \psi_i$ and $\bar{l}_{ip} \overset{a.s.}{\to} 1$ which implies that
     $$ \frac{l_{ip}}{\bar{l}_{ip}} \overset{a.s.}{\to} \psi_i \implies \left| \frac{l_{ip}}{\bar{l}_{ip}} - \psi_i \right| \overset{a.s.}{\to} 0. $$
     Hence the claim is proved.

 Next using the monotonicity of $F_c(.)$ and $\lambda_i$'s we have
 $$A^{**}_j-A^{**}_k \overset{a.s.}{\to} \sum_{i=j+1}^k F_c(\lambda_i) -  F_c(1+\sqrt{c}) > (k-j)[F_c(\lambda_k) -  F_c(1+\sqrt{c})]. $$
 We know that $\lambda_k > 1+\sqrt{c}$ so we have that $F_c(\lambda_k) -  F_c(1+\sqrt{c}) > 0$ hence $\mathbb{P}(\lim_{n\to \infty}(A^{**}_j - A^{**}_k )> 0) =1 $  $\forall j <k$ this implies that $\mathbb{P}(\lim _{n\to \infty}\hat{k}^{**}_A \geq k) = 1$ or a weaker condition that $\mathbb{P}(\hat{k}^{**}_A \geq k) \overset{n\to\infty}{\to} 1$. \\

\textbf{Case II:} When $j>k$ and $j$ is bounded, i.e. $j<s$,
\begin{align*}
     A^{**}_j-A^{**}_k &= \sum_{i=k+1}^j (A^{**}_{i} - A^{**}_{i-1}) \\
     \begin{split}
     &= \sum_{i=k+1}^j 
     -(p-i+1)\log\left\{1-\frac{1}{p-i+1} \left(1-\frac{l_{ip}}{\bar{l}_{ip}}\right)  \right\}\\
     & -\log\bar{l}_{ip} + \log l_{ip} + \frac{1}{p/n}F_{p/n}(1+\sqrt{p/n}+\delta_n)\frac{p-i+1}{n} 
    \end{split}\\
\end{align*}
Lets look at $h n^{2/3}(A^{**}_j-A^{**}_k)$ where $h \to 0$ as $n \to \infty$, pick the $ith$ term, we divide the expression into sum of four parts as follows
\begin{enumerate}
    \item  $-h n^{2/3} \log \bar{l}_{ip}$
    \item  $h n^{2/3} (\log l_{ip} - \log b) $
    \item  $h n^{2/3}\left[-(p-i+1)\log\left\{1-\frac{1}{p-i+1} \left(1-\frac{l_{ip}}{\bar{l}_{ip}}\right)  \right\} -(1-b)\right]$
    \item $h n^{2/3}\left(\frac{1}{p/n}F_{p/n}(1+\sqrt{p/n}+\delta_n)\frac{p-i+1}{n} + \log b +1-b\right)$
\end{enumerate}
Here $b:=\psi_c(1+\sqrt{c}) = (1+\sqrt{c})^2$.\\
Let us begin analysing each of these terms one by one 
\begin{enumerate}
    \item $-h n^{2/3} \log \bar{l}_{ip}$\\
    Using Lemma \ref{LiHanYao}: 
    $$trace(S_n) - trace(\Sigma_p) \to N(0,\sigma_0^2)$$
    $\implies$ $p (\frac{\sum_{j=1}^p l_{jp}}{p} - \frac{\sum_{j=1}^p \lambda_{j}}{p}) = O_p(1)$\\
    $\implies$ $p (\frac{\sum_{j=1}^i l_{jp} + \sum_{j=i+1}^p l_{jp} }{p} - \frac{\sum_{j=1}^i \lambda_{j} + (p-i)}{p}) = O_p(1)$\\
    $\implies$ $p\left(
    (\frac{p-i}{p})[\bar{l}_{ip}-1]+\frac{\sum_{j=1}^k l_{jp}-\sum_{j=1}^k \lambda_{j}}{p} + \frac{\sum_{j=k+1}^i (l_{jp}-\lambda_{j})}{p})\right) = O_p(1)$\\
    $\implies$ $(p-i)[\bar{l}_{ip}-1] + \sum_{j=1}^k (l_{jp}-\lambda_j) +\sum_{j=k+1}^i (l_{jp}-1) = O_p(1)$\\
    
    We know that $l_{jp} \to \psi_j=\psi_c(\lambda_j)$ $a.s.$ if $j\leq k$ otherwise  $l_{jp} \to b$ $a.s.$ when $j >k$ and $j$ is finite.
    $$\sum_{j=1}^k (l_{jp}-\lambda_j) +\sum_{j=k+1}^i (l_{jp}-1) \overset{a.s.}{\to} \sum_{j=1}^k (\psi_j -\lambda_j) + (i-k)(b-1). $$
    
    This along with the previous fact implies that
    $$(p-i)[\bar{l}_{ip}-1] = O_p(1).$$
    So $$h n^{2/3}[\bar{l}_{ip}-1] = \frac{h n^{2/3}}{p-i}(p-i)[\bar{l}_{ip}-1],$$
    coupled with the fact that $p/n \to c>0$ and Slutsky's theorem we have that 
    $$h n^{2/3}[\bar{l}_{ip}-1] \overset{p}{\to} 0.$$
    Next by mean value theorem we have that:
    $$\log x - \log 1 = \frac{1}{\zeta(x)}(x-1) \qquad \text{ where } \min(1,x)<\zeta(x)<\max(1,x).$$
    Using the above we can infer that 
    $$ h n^{2/3}\log \bar{l}_{ip} =  \frac{h n^{2/3}}{\zeta(\bar{l}_{ip})}(\bar{l}_{ip} - 1).$$
    As $\bar{l}_{ip} \overset{p}{\to} 1$. Using sandwich theorem we have:
     $$\zeta(\bar{l}_{ip}) \overset{p}{\to} 1. $$
     Therefore using Slutsky's Theorem we have that:
     $$-h n^{2/3}\log \bar{l}_{ip} = -\frac{h n^{2/3}}{\zeta(\bar{l}_{ip})}(\bar{l}_{ip} - 1) \overset{p}{\to} 0. $$
     
     \item $h n^{2/3} (\log l_{ip} - \log b)$\\
     By mean value theorem we have that:
     $$\log x - \log b = \frac{1}{\zeta(x)}(x-b) \qquad \text{ where } \min(b,x)<\zeta(x)<\max(b,x).$$
     Therefore
     $$\log l_{ip} - \log b= \frac{1}{\zeta(l_{ip})}(l_{ip}-b).$$
     Using Lemma \ref{passemierlemma}
     $$n^{2/3}(l_{ip}-b) =O_p(1) \quad \forall k <i \leq s.$$
     Also as $l_{ip} \overset{a.s.}{\to} b$ so using sandwich theorem we have $\zeta(l_{ip})\overset{a.s.}{\to} b$.
     Therefore we have that:
     $$h n^{2/3} (\log l_{ip} - \log b) = h n^{2/3}(l_{ip}-b)\frac{1}{\zeta(l_{ip})},$$
     so using the above fact we have that
     $$hn^{2/3}(l_{ip}-b) \overset{p}{\to}, 0$$
     which implies $$h n^{2/3} (\log l_{ip} - \log b) \overset{p}{\to} 0.$$

     \item $h n^{2/3}\left[-(p-i+1)\log\left\{1-\frac{1}{p-i+1} \left(1-\frac{l_{ip}}{\bar{l}_{ip}}\right)  \right\} -(1-b)\right]$\\
     
     Observe that,
     $h n^{2/3}\left|-(p-i+1)\log\left\{1-\frac{1}{p-i+1} \left(1-\frac{l_{ip}}{\bar{l}_{ip}}\right)  \right\} -(1-b)\right|$ \\
     $\leq$ $h n^{2/3}\left|-(p-i+1)\log\left\{1-\frac{1}{p-i+1} \left(1-\frac{l_{ip}}{\bar{l}_{ip}}\right)  \right\} - \left(1-\frac{l_{ip}}{\bar{l}_{ip}}\right)\right|$ +$\left| \frac{l_{ip}}{\bar{l}_{ip}} - b\right|.$\\
     
     Using the taylors series theorem we get under some conditions on $f$ we have
     $$f(x+h_0) = f(x) + f'(x)h_0 + \frac{f''(\zeta)}{2}h_0^2 \quad \zeta \in (\min(x,x+h_0),\max(x.x+h_0))).$$
     Choose 
     $$f(x) =\log(1-x)\qquad x=0  \qquad h_0 = \frac{1}{p-i+1}\left(1-\frac{l_{ip}}{\bar{l}_{ip}}\right)$$
     As in case I we finally have that
     $$\log\left\{1-\frac{1}{p-i+1} \left(1-\frac{l_{ip}}{\bar{l}_{ip}}\right)  \right\} = -\frac{1}{p-i+1}\left(1-\frac{l_{ip}}{\bar{l}_{ip}}\right)+ \frac{h_0^2}{2(1-\zeta)^2}$$
     where
     $\zeta \in (\min(0,h_0),\max(0,h_0))$\\
     So we have that 
     \begin{align*}
     \left|-(p-i+1)\log\left\{1-\frac{1}{p-i+1} \left(1-\frac{l_{ip}}{\bar{l}_{ip}}\right)  \right\} - \left(1-\frac{l_{ip}}{\bar{l}_{ip}}\right)\right| = \left|(p-i+1)\frac{h_0^2}{2(1-\zeta)^2}\right|
     \end{align*}
     Which is equal to
     $$\frac{1}{(1-\zeta)^2}\left(1-\frac{l_{ip}}{\bar{l}_{ip}}\right)^2\frac{1}{p-i+1}.$$
     Lets look at the R.H.S, we already know that $l_{ip}\overset{a.s.}{\to}b$ and $\bar{l}_{ip}\overset{a.s.}{\to}1.$ \\
     Using these we infer that
     $$\left(1-\frac{l_{ip}}{\bar{l}_{ip}}\right)^2 \overset{a.s}{\to} (1-b)^2 \quad \& \quad h = \frac{1}{p-i+1}\left(1-\frac{l_{ip}}{\bar{l}_{ip}}\right)\overset{a.s}{\to} 0.$$
     Now using sandwich theorem we have that $\zeta \overset{a.s}{\to} 0 $ as $\zeta \in (\min(0,h_0),\max(0,h_0))$ 
     so, we have that $\frac{1}{(1-\zeta)^2} \overset{a.s}{\to} 1$.\\
     Therefore,
     $$\frac{1}{(1-\zeta)^2}\left(1-\frac{l_{ip}}{\bar{l}_{ip}}\right)^2 \overset{a.s}{\to} (1-b)^2.$$
     So,
     $$h n^{2/3}\left|-(p-i+1)\log\left\{1-\frac{1}{p-i+1} \left(1-\frac{l_{ip}}{\bar{l}_{ip}}\right)  \right\} - \left(1-\frac{l_{ip}}{\bar{l}_{ip}}\right)\right| = \frac{\left(1-\frac{l_{ip}}{\bar{l}_{ip}}\right)^2}{(1-\zeta)^2}\frac{hn^{2/3}}{p-i+1}.$$
     As $\frac{hn^{2/3}}{p-i+1} \to 0$ coupled with previous fact, we can conclude that 
     $$h n^{2/3}\left|-(p-i+1)\log\left\{1-\frac{1}{p-i+1} \left(1-\frac{l_{ip}}{\bar{l}_{ip}}\right)  \right\} - \left(1-\frac{l_{ip}}{\bar{l}_{ip}}\right)\right| \overset{P}{\to} 0.$$
     Next let us deal with 
     \begin{align*}
         hn^{2/3}\left| \frac{l_{ip}}{\bar{l}_{ip}} - b\right| &= hn^{2/3}\left| \frac{l_{ip}}{\bar{l}_{ip}}-\frac{b}{\bar{l}_{ip}}+\frac{b}{\bar{l}_{ip}} - b\right|\\
         &\leq hn^{2/3}\left| \frac{l_{ip}}{\bar{l}_{ip}}-\frac{b}{\bar{l}_{ip}}\right|+hn^{2/3}\left|\frac{b}{\bar{l}_{ip}} - b\right|\\
         & \leq hn^{2/3}\frac{|l_{ip}-b|}{|\bar{l}_{ip}|} +  hn^{2/3}b\frac{|\bar{l}_{ip}-1|}{|\bar{l}_{ip}|}
     \end{align*}
     Using Lemma \ref{passemierlemma}
     $$n^{2/3}(l_{ip}-b) = O_p(1) \quad \forall k <i \leq s$$
     Therefore
     $$hn^{2/3}(l_{ip}-b) \overset{P}{\to} 0 \quad \&\quad \bar{l}_{ip} \overset{a.s.}{\to} 1 \implies hn^{2/3}\frac{|l_{ip}-b|}{|\bar{l}_{ip}|} \overset{P}{\to} 0.$$
     We have already observed that 
     $$h n^{2/3}[\bar{l}_{ip}-1] \overset{p}{\to} 0$$
     Therefore
     $$h n^{2/3}[\bar{l}_{ip}-1] \overset{p}{\to} 0 \quad \&\quad \bar{l}_{ip} \overset{a.s.}{\to} 1 \implies b hn^{2/3}\frac{|l_{ip}-1|}{|\bar{l}_{ip}|} \overset{P}{\to} 0.$$
     Hence
     $$hn^{2/3}\left| \frac{l_{ip}}{\bar{l}_{ip}} - b\right| \overset{P}{\to} 0.$$
     So finally we have that
     $$h n^{2/3}\left[-(p-i+1)\log\left\{1-\frac{1}{p-i+1} \left(1-\frac{l_{ip}}{\bar{l}_{ip}}\right)  \right\} -(1-b)\right]\overset{P}{\to} 0.$$
     
     \item Looking at the last term 
     $$h n^{2/3}\left(\frac{1}{p/n}F_{p/n}(1+\sqrt{p/n}+\delta_n)\frac{p-i+1}{n} + \log b +1-b\right)$$
     Call
     $$\alpha(p,n) = \frac{1}{p/n}F_{p/n}(1+\sqrt{p/n}+\delta_n) \quad \& \quad \alpha(c) =\frac{1}{c} F_c(\bbp)$$
     therefore the last term can be written as
     \begin{align*}
         h n^{2/3}\left(\alpha(p,n)\frac{p}{n} + \log b +1-b\right) - h n^{2/3}\alpha(p,n)\frac{i-1}{n}.
     \end{align*}
     Observe that
     $$h n^{2/3}\frac{i-1}{n} \to 0\quad \&\quad \alpha(p,n)\to\alpha(c) \implies h n^{2/3}\alpha(p,n)\frac{i-1}{n} \to 0$$
     Therefore asymptotically only the first term matters, which is
     $$h n^{2/3}\left(\alpha(p,n)\frac{p}{n} + \log b +1-b\right) = h n^{2/3}\left(F_{p/n}(1+\sqrt{p/n}+\delta_n)  -F_c(1+\sqrt{c})\right).$$
     Now $\delta_n$ is chosen such that this term above goes to $\infty$. 
\end{enumerate}
As a result we have that
     $$hn^{2/3}(A^{**}_j-A^{**}_k) \overset{p}{\to} \infty \,  \implies \mathbb{P}(A^{**}_j>A^{**}_k) \to 1 \quad s\geq j>k$$
    \\
Now combining both the cases we have $$\mathbb{P}(A^{**}_j > A^{**}_k) \to 1 \quad \forall j\neq k \quad j\leq s$$\\
As there are finitely many $j$'s this implies that
$$\mathbb{P}(A^{**}_j >A^{**}_k \,\,\forall j\neq k \,\, j\leq s) \to 1$$
Hence $\mathbb{P}(\hat{k}^{*}_A = k) \to 1 \implies \hat{k}^{*}_A \overset{P}{\to} k$.

\phantomsection
\label{weakconsproof2}
\textbf{Proof of the Theorem \ref{weakconsistencythm2}:} 
We will show that $\mathbb{P}(\hat{k}^{**}_{\tilde{A}} = k)\to 1$ as $n \to \infty$. We will show that by breaking the problem into two cases as before $j<k$ and $j>k$, in both these cases we will show that $\mathbb{P}(\tilde{A}^{**}_j >\tilde{A}^{**}_k)\to 1$ as $n \to \infty$.\\ 
 \textbf{Case I:} When $j<k$,
 \begin{align*}
     \tilde{A}^{**}_j-\tilde{A}^{**}_k &= \sum_{i=j+1}^k (\tilde{A}^{**}_{i-1} - \tilde{A}^{**}_i) \\
     \begin{split}
     &= \sum_{i=j+1}^k 
     (n-i)\log\left\{1-\frac{1}{n-i} \left(1-\frac{l_{ip}}{\bar{l}_{ip}}\right)  \right\}\\
     & +\log\bar{l}_{ip} - \log l_{ip} - Q_{p/n}(1+\sqrt{p/n}+\delta_n)\frac{n-i}{p} 
    \end{split}\\
    &\overset{a.s.}{\to} \sum_{i=j+1}^k \psi_i/c - 1 - \log(\psi_i/c) - \frac{1}{c}Q_c(\bbp)\\
    &= \sum_{i=j+1}^k\frac{1}{c} Q_c(\lambda_i) - \frac{1}{c}Q_c(\bbp) 
 \end{align*}
 Only step that needs to be justified is the third step , specifically\\
 \textbf{Claim:}
 $$(n-i)\log\left\{1-\frac{1}{n-i} \left(1-\frac{l_{ip}}{\bar{l}_{ip}}\right)\right\} \overset{a.s.}{\to} \psi_i/c -1.$$
Observe that,
     $\left|(n-i)\log\left\{1-\frac{1}{n-i} \left(1-\frac{l_{ip}}{\bar{l}_{ip}}\right)  \right\} -(\psi_i/c-1)\right|$ \\
     $\leq$ $\left|-(n-i)\log\left\{1-\frac{1}{n-i} \left(1-\frac{l_{ip}}{\bar{l}_{ip}}\right)  \right\} - \left(1-\frac{l_{ip}}{\bar{l}_{ip}}\right)\right|$ +$\left| \frac{l_{ip}}{\bar{l}_{ip}} - \frac{\psi_i}{c} \right|$\\
     Using Taylor's Series theorem with
     $$f(x) =\log(1-x)\qquad x=0  \qquad h_0 = \frac{1}{n-i}\left(1-\frac{l_{ip}}{\bar{l}_{ip}}\right)$$
     We have that
     $$\log\left\{1-\frac{1}{n-i} \left(1-\frac{l_{ip}}{\bar{l}_{ip}}\right)  \right\} = -\frac{1}{n-i}\left(1-\frac{l_{ip}}{\bar{l}_{ip}}\right)+ \frac{h_0^2}{2(1-\zeta)^2},$$
     where
     $\zeta \in (\min(0,h_0),\max(0,h_0)).$\\
     So we have
     $$-(n-i)\log\left\{1-\frac{1}{n-i} \left(1-\frac{l_{ip}}{\bar{l}_{ip}}\right)  \right\} = \left(1-\frac{l_{ip}}{\bar{l}_{ip}}\right)-(n-i) \frac{h_0^2}{2(1-\zeta)^2}.$$
     Hence
     \begin{align*}
     \left|-(n-i)\log\left\{1-\frac{1}{n-i} \left(1-\frac{l_{ip}}{\bar{l}_{ip}}\right)  \right\} - \left(1-\frac{l_{ip}}{\bar{l}_{ip}}\right)\right| &= \left|(n-i)\frac{h_0^2}{2(1-\zeta)^2}\right|\\
     =& \frac{1}{(1-\zeta)^2}\left(1-\frac{l_{ip}}{\bar{l}_{ip}}\right)^2\frac{1}{n-i}
     \end{align*}
     Lets look at the R.H.S, we know from Lemma \ref{asconvergence} that $l_{ip}\overset{a.s.}{\to}\psi_i$ .\\
     Next we look at 
     $$\bar{l}_{ip} = \frac{1}{n-1-i}\sum_{t=i+1}^{n-1} l_{tp}$$
     We already know by M-P law that 
     $$\frac{1}{p-i}\sum_{t=i+1}^{p} l_{tp} =\frac{n-1-i}{p-i} \frac{1}{n-1-i}\sum_{t=i+1}^{n-1} l_{tp} \overset{a.s.}{\to} \mu_{MP} = 1. $$
     where $\mu_{MP}$ is the mean of the Marcenko-Pastur distribution.
     Therefore
     $$\frac{n-1-i}{p-i} \frac{1}{n-1-i}\sum_{t=i+1}^{n-1} l_{tp} = \frac{n-1-i}{p-i} \bar{l}_{ip} \overset{a.s.}{\to} 1 \implies \bar{l}_{ip} \overset{a.s.}{\to} c. $$
     Using these we infer that
     $$\left(1-\frac{l_{ip}}{\bar{l}_{ip}}\right)^2 \overset{a.s}{\to} \left(1-\frac{\psi_i}{c}\right)^2 \quad \& \quad h_0 = \frac{1}{n-i}\left(1-\frac{l_{ip}}{\bar{l}_{ip}}\right)\overset{a.s}{\to} 0. $$
     Now using sandwich theorem we have that $\zeta \overset{a.s}{\to} 0 $ as $\zeta \in (\min(0,h_0),\max(0,h_0))$ 
     so, we have that $\frac{1}{(1-\zeta)^2} \overset{a.s}{\to} 1$.\\
     Therefore,
     $$\frac{1}{(1-\zeta)^2}\left(1-\frac{l_{ip}}{\bar{l}_{ip}}\right)^2 \overset{a.s}{\to} \left(1-\frac{\psi_i}{c}\right)^2.$$
     So,
     $$\left|-(n-i)\log\left\{1-\frac{1}{n-i} \left(1-\frac{l_{ip}}{\bar{l}_{ip}}\right)  \right\} - \left(1-\frac{l_{ip}}{\bar{l}_{ip}}\right)\right| = \frac{\left(1-\frac{l_{ip}}{\bar{l}_{ip}}\right)^2}{(1-\zeta)^2}\frac{1}{n-i}$$
     As $\frac{1}{n-i} \to 0$ coupled with previous fact, we can conclude that 
     $$\left|-(n-i)\log\left\{1-\frac{1}{n-i} \left(1-\frac{l_{ip}}{\bar{l}_{ip}}\right)  \right\} - \left(1-\frac{l_{ip}}{\bar{l}_{ip}}\right)\right| \overset{a.s.}{\to} 0,$$
     As already observed $l_{ip} \overset{a.s.}{\to} \psi_i$ and $\bar{l}_{ip} \overset{a.s.}{\to} c$ which implies that
     $$ \frac{l_{ip}}{\bar{l}_{ip}} \overset{a.s.}{\to} \frac{\psi_i}{c} \implies \left| \frac{l_{ip}}{\bar{l}_{ip}} - \frac{\psi_i}{c} \right| \overset{a.s.}{\to} 0.$$ 
     Hence the claim is proved.

 Next using the monotonicity of $Q_c(.)$ and $\lambda_i$'s
 $$\tilde{A}^{**}_j-\tilde{A}^{**}_k \overset{a.s.}{\to} \sum_{i=j+1}^k \frac{1}{c}(Q_c(\lambda_i) -  Q_c(\bbp)) > (k-j)\frac{1}{c}[Q_c(\lambda_k) -  Q_c(\bbp]. $$
 We know that $\lambda_k > 1+\sqrt{c}$ so we have that $Q_c(\lambda_k) -  Q_c(1+\sqrt{c}) > 0$ hence $\mathbb{P}(\lim_{n\to \infty}(\tilde{A}^{**}_j - \tilde{A}^{**}_k )> 0) =1 $  $\forall j <k$ this implies that $\mathbb{P}(\lim _{n\to \infty}\hat{k}^{**}_{\tilde{A}} \geq k) = 1$ or a weaker condition that $\mathbb{P}(\hat{k}^{**}_{\tilde{A}} \geq k) \overset{n\to\infty}{\to} 1$. \\

\textbf{Case II:} When $j>k$ and $j$ is bounded, i.e. $j<s$,
\begin{align*}
     \tilde{A}^{**}_j-\tilde{A}^{**}_k &= \sum_{i=k+1}^j (\tilde{A}^{**}_{i} - \tilde{A}^{**}_{i-1}) \\
     \begin{split}
     &= \sum_{i=k+1}^j 
     -(n-i)\log\left\{1-\frac{1}{n-i} \left(1-\frac{l_{ip}}{\bar{l}_{ip}}\right)  \right\}\\
     & -\log\bar{l}_{ip} + \log l_{ip} + Q_{p/n}(1+\sqrt{p/n}+\delta_n)\frac{n-i}{p} 
    \end{split}\\
\end{align*}
Lets look at $h n^{2/3}(\tilde{A}^{**}_j-\tilde{A}^{**}_k)$ where $h \to 0$ as $n \to \infty$, pick the $ith$ term, we divide the expression into sum of four parts as follows
\begin{enumerate}
    \item  $-h n^{2/3} (\log \bar{l}_{ip} - \log c)$
    \item  $h n^{2/3} (\log l_{ip} - \log b) $
    \item  $h n^{2/3}\left[-(n-i)\log\left\{1-\frac{1}{n-i} \left(1-\frac{l_{ip}}{\bar{l}_{ip}}\right)  \right\} -(1-\frac{b}{c})\right]$
    \item $h n^{2/3}\left(Q_{p/n}(1+\sqrt{p/n}+\delta_n)\frac{n-i}{p} + \log \frac{b}{c} +1-\frac{b}{c}\right)$
\end{enumerate}
Here $b:=\psi_c(1+\sqrt{c}) = (1+\sqrt{c})^2$.\\
Let us begin analysing each of these terms one by one 
\begin{enumerate}
    \item $-h n^{2/3}( \log \bar{l}_{ip} -\log c)$\\
    
    From the proof of Theorem \ref{weakconsitencythm1}
    we have already shown that
    $$hn^{2/3}\left[\frac{1}{p-i}\sum_{t=i+1}^p l_{tp}-1\right] = o_p(1)$$
    Using which we infer that
    $$hn^{2/3} \left( \frac{n-i-1}{p-i}\bar{l}_{ip}- 1 \right)= o_p(1)$$
    As $\frac{p-i}{n-i-1} \to c $ , we have that:
    $$\frac{p-i}{n-i-1}hn^{2/3} \left( \frac{n-i-1}{p-i}\bar{l}_{ip}- 1 \right)= o_p(1)$$
    This implies
    $$hn^{2/3}\left(\bar{l}_{ip}- \frac{p-i}{n-i-1} \right)= o_p(1),$$
    i.e. 
    $$hn^{2/3} (\bar{l}_{ip}- c)+hn^{2/3}\left(c-\frac{p-i}{n-i-1} \right) = o_p(1)$$
    We have that $$p/n = c +O\left(n^{-2/3}\right) \implies hn^{2/3}\left(c-\frac{p-i}{n-i-1} \right) = o(1) $$
    As a result
    $$hn^{2/3} (\bar{l}_{ip}- c)=o_p(1)$$
    
    Next by mean value theorem we have that:
    $$\log x - \log c = \frac{1}{\zeta(x)}(x-c) \qquad \text{ where } \min(c,x)<\zeta(x)<\max(c,x)$$
    Using the above we can infer that 
    $$ h n^{2/3}(\log \bar{l}_{ip}-\log c)  =  \frac{h n^{2/3}}{\zeta(\bar{l}_{ip})}(\bar{l}_{ip} - c)$$
    As already observed during the proof of Theorem \ref{weakconsitencythm1} $\bar{l}_{ip} \overset{p}{\to} c$. Using sandwich theorem we have:
     $$\zeta(\bar{l}_{ip}) \overset{p}{\to} c $$
     Therefore using Slutsky's Theorem we have that:
     $$-h n^{2/3}(\log \bar{l}_{ip}-\log c) = -\frac{h n^{2/3}}{\zeta(\bar{l}_{ip})}(\bar{l}_{ip} - c) \overset{p}{\to} 0 $$
    
    \item $h n^{2/3} (\log l_{ip} - \log b)$\\
     The proof of this part is exactly the same as that of Theorem \ref{weakconsitencythm1}, where we showed that
     $$h n^{2/3} (\log l_{ip} - \log b) \overset{p}{\to} 0$$
     
     \item $h n^{2/3}\left[-(n-i)\log\left\{1-\frac{1}{n-i} \left(1-\frac{l_{ip}}{\bar{l}_{ip}}\right)  \right\} -(1-\frac{b}{c})\right]$\\
     
     Observe that,
     $h n^{2/3}\left|-(n-i)\log\left\{1-\frac{1}{n-i} \left(1-\frac{l_{ip}}{\bar{l}_{ip}}\right)  \right\} -(1-\frac{b}{c})\right|$ \\
     $\leq$ $h n^{2/3}\left|-(n-i)\log\left\{1-\frac{1}{n-i} \left(1-\frac{l_{ip}}{\bar{l}_{ip}}\right)  \right\} - \left(1-\frac{l_{ip}}{\bar{l}_{ip}}\right)\right|$ +$\left| \frac{l_{ip}}{\bar{l}_{ip}} - \frac{b}{c}\right|$\\
     
     Using the taylors series theorem we get under some conditions on $f$ we have
     $$f(x+h_0) = f(x) + f'(x)h_0 + \frac{f''(\zeta)}{2}h_0^2 \quad \zeta \in (\min(x,x+h_0),\max(x.x+h_0)))$$
     Choose 
     $$f(x) =\log(1-x)\qquad x=0  \qquad h_0 = \frac{1}{n-i}\left(1-\frac{l_{ip}}{\bar{l}_{ip}}\right)$$
     As in case I we finally have that
     $$\log\left\{1-\frac{1}{n-i} \left(1-\frac{l_{ip}}{\bar{l}_{ip}}\right)  \right\} = -\frac{1}{n-i}\left(1-\frac{l_{ip}}{\bar{l}_{ip}}\right)+ \frac{h_0^2}{2(1-\zeta)^2}$$
     where
     $\zeta \in (\min(0,h_0),\max(0,h_0))$\\
     So we have that 
     \begin{align*}
     \left|-(n-i)\log\left\{1-\frac{1}{n-i} \left(1-\frac{l_{ip}}{\bar{l}_{ip}}\right)  \right\} - \left(1-\frac{l_{ip}}{\bar{l}_{ip}}\right)\right| &= \left|(n-i)\frac{h_0^2}{2(1-\zeta)^2}\right|\\
     =& \frac{1}{(1-\zeta)^2}\left(1-\frac{l_{ip}}{\bar{l}_{ip}}\right)^2\frac{1}{n-i}
     \end{align*}
     Lets look at the R.H.S, we already know that $l_{ip}\overset{a.s.}{\to}b$ and $\bar{l}_{ip}\overset{a.s.}{\to}c$ \\
     Using these we infer that
     $$\left(1-\frac{l_{ip}}{\bar{l}_{ip}}\right)^2 \overset{a.s}{\to} (1-b/c)^2 \quad \& \quad h = \frac{1}{n-i}\left(1-\frac{l_{ip}}{\bar{l}_{ip}}\right)\overset{a.s}{\to} 0 $$
     Now using sandwich theorem we have that $\zeta \overset{a.s}{\to} 0 $ as $\zeta \in (\min(0,h_0),\max(0,h_0))$ 
     so, we have that $\frac{1}{(1-\zeta)^2} \overset{a.s}{\to} 1$.\\
     Therefore,
     $$\frac{1}{(1-\zeta)^2}\left(1-\frac{l_{ip}}{\bar{l}_{ip}}\right)^2 \overset{a.s}{\to} \left(1-\frac{b}{c} \right)^2$$
     So,
     $$h n^{2/3}\left|-(n-i)\log\left\{1-\frac{1}{n-i} \left(1-\frac{l_{ip}}{\bar{l}_{ip}}\right)  \right\} - \left(1-\frac{l_{ip}}{\bar{l}_{ip}}\right)\right| = \frac{\left(1-\frac{l_{ip}}{\bar{l}_{ip}}\right)^2}{(1-\zeta)^2}\frac{hn^{2/3}}{n-i}$$
     As $\frac{hn^{2/3}}{n-i} \to 0$ coupled with previous fact, we can conclude that 
     $$h n^{2/3}\left|-(n-i)\log\left\{1-\frac{1}{n-i} \left(1-\frac{l_{ip}}{\bar{l}_{ip}}\right)  \right\} - \left(1-\frac{l_{ip}}{\bar{l}_{ip}}\right)\right| \overset{P}{\to} 0,$$
     Next let us deal with 
     \begin{align*}
         hn^{2/3}\left| \frac{l_{ip}}{\bar{l}_{ip}} - \frac{b}{c}\right| &= hn^{2/3}\left| \frac{l_{ip}}{\bar{l}_{ip}}-\frac{b}{\bar{l}_{ip}}+\frac{b}{\bar{l}_{ip}} - \frac{b}{c}\right|\\
         &\leq hn^{2/3}\left| \frac{l_{ip}}{\bar{l}_{ip}}-\frac{b}{\bar{l}_{ip}}\right|+hn^{2/3}\left|\frac{b}{\bar{l}_{ip}} - \frac{b}{c}\right|\\
         & \leq hn^{2/3}\frac{|l_{ip}-b|}{|\bar{l}_{ip}|} +  hn^{2/3}\frac{b}{c}\frac{|\bar{l}_{ip}-c|}{|\bar{l}_{ip}|}
     \end{align*}
     Using Lemma \ref{passemierlemma}
     $$n^{2/3}(l_{ip}-b) = O_p(1) \quad \forall k <i \leq s$$
     Therefore
     $$hn^{2/3}(l_{ip}-b) \overset{P}{\to} 0 \quad \&\quad \bar{l}_{ip} \overset{a.s.}{\to} c \implies hn^{2/3}\frac{|l_{ip}-b|}{|\bar{l}_{ip}|} \overset{P}{\to} 0$$
     We have already observed that 
     $$h n^{2/3}[\bar{l}_{ip}-1] \overset{p}{\to} 0$$
     Therefore
     $$h n^{2/3}[\bar{l}_{ip}-c] \overset{p}{\to} 0 \quad \&\quad \bar{l}_{ip} \overset{a.s.}{\to} c \implies  hn^{2/3}\frac{b}{c}\frac{|l_{ip}-c|}{|\bar{l}_{ip}|} \overset{P}{\to} 0$$
     Hence
     $$hn^{2/3}\left| \frac{l_{ip}}{\bar{l}_{ip}} - \frac{b}{c}\right| \overset{P}{\to} 0$$
     So finally we have that
     $$h n^{2/3}\left[-(p-i+1)\log\left\{1-\frac{1}{p-i+1} \left(1-\frac{l_{ip}}{\bar{l}_{ip}}\right)  \right\} -\left(1-\frac{b}{c}\right)\right]\overset{P}{\to} 0$$

     \item Looking at the last term 
     $$h n^{2/3}\left(Q_{p/n}(1+\sqrt{p/n}+\delta_n)\frac{n-i}{p} + \log \frac{b}{c} +1-\frac{b}{c}\right)$$
     Call
     $$\beta(p,n) = Q_{p/n}(1+\sqrt{p/n}+\delta_n) \quad \& \quad \beta(c) = Q_c(\bbp)$$
     therefore the last term can be written as
     \begin{align*}
         h n^{2/3}\left(\beta(p,n)\frac{n}{p} + \log \frac{b}{c} +1-\frac{b}{c}\right) - h n^{2/3}\beta(p,n)\frac{i}{p}
     \end{align*}
     Observe that
     $$h n^{2/3}\frac{i}{p} \to 0\quad \&\quad \beta(p,n)\to\beta(c) \implies  h n^{2/3}\beta(p,n)\frac{i}{p} \to 0$$
     Therefore asymptotically only the first term matters, which is
     \begin{align*}
         h n^{2/3}&\left(\beta(p,n)\frac{n}{p} + \log \frac{b}{c} +1-\frac{b}{c}\right) =\\ 
         &h n^{2/3}\left(\frac{1}{p/n}Q_{p/n}(1+\sqrt{p/n}+\delta_n)  -\frac{1}{c}Q_c(1+\sqrt{c})\right) \to \infty 
     \end{align*}
     as $n \to \infty$    
\end{enumerate}
As a result we have that
     $$hn^{2/3}(\tilde{A}^{**}_j-\tilde{A}^{**}_k) \overset{p}{\to} \infty \,  \implies \mathbb{P}(\tilde{A}^{**}_j>\tilde{A}^{**}_k) \to 1 \quad s\geq j>k$$
    \\
Now combining both the cases we have $$\mathbb{P}(\tilde{A}^{**}_j > \tilde{A}^{**}_k) \to 1 \quad \forall j\neq k \quad j\leq s$$\\
As there are finitely many $j$'s this implies that
$$\mathbb{P}(\tilde{A}^{**}_j >\tilde{A}^{**}_k \,\,\forall j\neq k \,\, j\leq s) \to 1$$
Hence $\mathbb{P}(\hat{k}^{*}_{\tilde{A}} = k) \to 1 \implies \hat{k}^{*}_{\tilde{A}} \overset{P}{\to} k$.

\subsection{Strong Consistency Proof}
\label{strongconsproof1}
\textbf{Proof of Theorem \ref{strongconsistencythm1}} (Strong Consistency): Here again as for weak consistency we look for two cases $j<k$ and $j>k$ and compare $A_j$ and $A_k$. Proof of Case I, i.e. $j<k$ is along the similar lines as in Theorem \ref{weakconsitencythm1} (weak consistency proof) where we have shown that 
\[
\mathbb{P}(\lim_{n\to \infty}\hat{k}^*_A \geq k ) = 1.
\]
\textbf{Case I:} When $j<k$,

 \begin{align*}
     A_j-A_k &= \sum_{i=j+1}^k (A_{i-1} - A_i) \\
     \begin{split}
     &= \sum_{i=j+1}^k 
     (p-i+1)\log\left\{1-\frac{1}{p-i+1} \left(1-\frac{l_{ip}}{\bar{l}_{ip}}\right)  \right\}\\
     & +\log\bar{l}_{ip} - \log l_{ip} - \frac{1}{p/n}F_{p/n}(1+\sqrt{p/n}+\delta)\frac{p-i+1}{n} 
    \end{split}\\
    &\overset{a.s.}{\to} \sum_{i=j+1}^k \psi_i - 1 - \log\psi_i - F_c(1+\sqrt{c}+\delta) \qquad \\
    &= \sum_{i=j+1}^k F_c(\lambda_i) - F_c(1+\sqrt{c}+\delta) 
 \end{align*}
 Using the monotonicity of $F_c(.)$ and $\lambda_i$'s
 $$A_j-A_k \overset{a.s.}{\to} \sum_{i=j+1}^k F_c(\lambda_i) - F_c(1+\sqrt{c}+\delta) > (k-j)[F_c(\lambda_k) - F_c(1+\sqrt{c}+\delta)]. $$
 Now suppose $\lambda_k < 1+\sqrt{c}+\delta$ then we have that $F_c(\lambda_k) - F_c(1+\sqrt{c}+\delta) <0$. So we have that
 $$A_{k-1} - A_k \overset{a.s.}{\to} F_c(\lambda_k) - F_c(1+\sqrt{c}+\delta) < 0,$$
 i.e. $A_{k-1}>A_k$ almost surely implying that $\hat{k}^{**}_A$ is not consistent.
 
On the other hand if we have that $\lambda_k > 1+\sqrt{c}+\delta$ so we have that $F_c(\lambda_k) - F_c(1+\sqrt{c}+\delta) > 0$ hence $\mathbb{P}(\lim_{n\to \infty}(A_j - A_k)> 0) =1 $  $\forall j <k$ this implies that $\mathbb{P}(\lim _{n\to \infty}\hat{k}^{**}_A \geq k) = 1$.
 
 \textbf{Case II:} When $j>k$ ,
\begin{align*}
     A_j-A_k &= \sum_{i=k+1}^j (A_{i} - A_{i-1}) \\
     \begin{split}
     &= \sum_{i=k+1}^j 
     -(p-i+1)\log\left\{1-\frac{1}{p-i+1} \left(1-\frac{l_{ip}}{\bar{l}_{ip}}\right)  \right\}\\
     & -\log\bar{l}_{ip} + \log l_{ip} + \frac{1}{p/n}F_{p/n}(1+\sqrt{p/n}+\delta)\frac{p-i+1}{n} 
    \end{split}\\
    &\sim \sum_{i=k+1}^j \left(1-\frac{l_{ip}}{\bar{l}_{ip}}\right) +\log\left(\frac{l_{ip}}{\bar{l}_{ip}}\right) + F_c(1+\sqrt{c}+\delta)\left(1-\frac{i}{p}\right)
\end{align*}
For$ k<j\leq j$, $l_{jp}\leq l_{ip} \leq l_{k+1,p}$. From Lemma \ref{asconvergence}, we have that $l_{k+1,p}$ and $l_{jp} \overset{a.s.}{\to} \mu_1 = b$ as $n \to \infty$ so  $l_{ip}\overset{a.s.}{\to} b $ it implies that a.s.
\begin{align*}
A_j -A_k &\sim (j-k)(1- b+ \log b + F_c(1+\sqrt{c}+\delta))\\
&=(j-k)(F_c(1+\sqrt{c}+\delta)-F_c(1+\sqrt{c}))>0.
\end{align*}
The last line is true because $ b-1- \log b = F_c(1+\sqrt{c})$ and because $F_c(.)$ is a monotonically increasing function.
Therefore we have that $\mathbb{P}(\lim_{n\to \infty}(A_j - A_k)> 0) =1 $  $\forall q>j>k$ this implies that $\mathbb{P}(\lim _{n\to \infty}\hat{k}^{**}_A \leq k) = 1$.\\
Combining Case I and II we have that $\hat{k}^{**}_A \overset{a.s.}{\to} k$. 
\end{document}